\documentclass[review,hidelinks,onefignum,onetabnum,letterpaper]{siamart220329}

\usepackage{graphicx}
\usepackage[pass,letterpaper]{geometry}
\usepackage{verbatim}
\usepackage{slashed}
\usepackage{graphicx}
\usepackage{subfigure}
\usepackage{amssymb}
\usepackage{subfigure}
\usepackage{color}
\usepackage{mathrsfs}
\usepackage{multirow}
\usepackage{listings}
\usepackage{array}
\usepackage{float}
\usepackage{framed}
\usepackage{appendix}
\newsiamthm{example}{Example}
\newsiamthm{remark}{Remark}

\allowdisplaybreaks

\newcommand{\PA}[1]{\partial #1}

\usepackage{algorithm}
\usepackage{algorithmic}
\usepackage{clrscode3e}
\usepackage{tabularx}
\usepackage{multirow}
\renewcommand{\vec}[1]{\boldsymbol{#1}}

\title{A volume-conserving ALE continuous Galerkin method with residual-based Dirichlet--Neumann fluxes for nonlinear free-surface waves}
\headers{A volume-conserving ALE-CG method}{Wang, Xia, and Xu}

\author{Fan Wang\thanks{School of Mathematical Sciences, University of Science
and Technology of China, Hefei, Anhui 230026, P.R. China. Email: ewangwf@mail.ustc.edu.cn.}
\and
Yinhua Xia\thanks{School of Mathematical Sciences, University of Science and Technology of China, Hefei, Anhui 230026, P.R. China.  E-mail: yhxia@ustc.edu.cn. The research of Yinhua Xia was partially supported by Anhui Provincial Natural Science Foundation grant No. 2408085J004  and NSFC grant No. 12671438. }
\and
Yan Xu\thanks{Corresponding author. School of Mathematical Sciences, University of Science and Technology of China, Hefei, Anhui 230026, P.R. China and Laoshan Laboratory, Qingdao 266237, P.R. China.  E-mail: yxu@ustc.edu.cn. The work of Yan Xu was partially supported by NSFC grant No. 12471347,  and Laoshan Laboratory (No. LSKJ202300305). } 
}
\date{}

\begin{document}
\nolinenumbers

\maketitle

\begin{abstract}
We develop an arbitrary Lagrangian--Eulerian continuous Galerkin finite
element method for two-dimensional fully nonlinear free-surface
potential flow. An isoparametric finite element method computes the
interior potential on the current fluid domain, while a residual-based
Dirichlet--Neumann flux drives the surface evolution. The residual-based flux satisfies a global discrete balance inherited from the interior variational formulation. When combined with
compatible mesh motion, surface projection, and quadrature, this balance
yields conservation of the discrete fluid volume at the semi-discrete
level. An optional volume-conserving sixth-order surface dissipation controls high-frequency
oscillations, and the surface equations are advanced with fourth-order
Runge--Kutta integration, updating the mesh and interior solution at every
stage. Numerical experiments confirm the expected spatial convergence for
linear free-surface waves and demonstrate reduced mesh sensitivity for
nonlinear waves over flat and curved bottoms. A rotating-wavemaker
example further illustrates the treatment of moving solid boundaries.
\end{abstract}



\section{Introduction}
Free-surface gravity waves arise in coastal and offshore engineering, sloshing, and laboratory wave generation. For smooth, non-breaking flows with weak viscous and rotational effects, the potential-flow model provides an effective description. Its numerical solution couples an elliptic problem in the fluid interior with nonlinear evolution equations on an unknown moving surface. A consistent treatment of the potential, surface geometry, and boundary motion is therefore essential for accurate wave simulation.

Zakharov's surface-variable formulation \cite{Zakharov1968} expresses the evolution in terms of the surface elevation and the trace of the velocity potential. The Dirichlet--Neumann operator (DNO) represents the harmonic extension into the fluid interior. Craig and Sulem \cite{CraigSulem1993} developed a spectral method based on an expansion of this operator. High-order spectral methods were also developed by Dommermuth and Yue \cite{DommermuthYue1987} and West et al.\ \cite{WestEtAl1987}. These methods exploit surface representations to compute smooth nonlinear wave fields accurately and efficiently. Extending them beyond the underlying reference geometry, however, requires additional treatment of the boundary configuration. The conditioning of high-order perturbative DNO evaluations also requires careful consideration, as analyzed by Nicholls and Reitich \cite{NichollsReitich2001}.

Discretizations in the physical domain provide another route to representing the fluid boundaries. Boundary integral methods have been widely used for nonlinear water waves \cite{GrilliEtAl1989}; a recent arbitrary Lagrangian--Eulerian (ALE)  regularized boundary-integral method addresses nonlinear waves over complex topography and wave--structure interaction \cite{TsaoKees2023}. Finite element methods combine geometric flexibility with a variational description of the interior problem. Cai et al.\ \cite{CaiEtAl1998} developed a finite element method for fully nonlinear potential-flow waves using a transformation of the fluid domain, while Wu and Eatock Taylor \cite{wu1994finite} developed finite element formulations for two-dimensional nonlinear transient waves, including wavemaker and rectangular-basin problems. Gagarina et al.\ \cite{GagarinaEtAl2014} coupled free-surface evolution and mesh deformation through a variational space--time formulation with symplectic time integration. Brink et al.\ \cite{BrinkEtAl2017} established a Hamiltonian finite element discretization on general deforming meshes, including prescribed boundary motion. More recent finite element and high-order formulations include a two-dimensional $\sigma$-transform method \cite{YangZhuRen2024} and a hybrid-spectral method for fully nonlinear potential flow \cite{MelanderEngsigKarup2024}. Discontinuous Galerkin formulations have also been developed for linear free-surface waves \cite{van2005discontinuous,tomar2007runge}. These studies demonstrate both the geometric flexibility of physical-domain discretizations and the importance of retaining the structure of the water-wave equations.

The ALE framework, reviewed by Donea et al.\ \cite{DoneaEtAl2004} and supported by the finite element foundations of Hughes et al.\ \cite{HughesEtAl1981}, permits mesh motion to be chosen independently of fluid particle motion. For free-surface flows, this freedom accommodates surface deformation and prescribed wall motion within a conforming mesh. The quasi-ALE finite element method of Ma and Yan \cite{MaYan2006} likewise uses a fixed initial mesh that is subsequently moved to conform to the evolving free surface and structures. Geometric flexibility alone, however, does not guarantee a discrete fluid-volume balance. Such a balance depends on the compatibility among the interior flux, surface evolution, and moving geometry. When numerical dissipation is introduced to control unresolved oscillations, its contribution to the volume rate must also be considered. These observations motivate a formulation in which both the flux balance and the effect of stabilization can be identified explicitly.

In this work, we develop an ALE continuous Galerkin finite element method for two-dimensional fully nonlinear free-surface gravity waves. A high-order isoparametric discretization solves the Laplace problem on the instantaneous fluid domain. The action of the DNO is then recovered from the finite element residual at the free-surface degrees of freedom. This construction uses the same variational problem for both the interior potential and the surface flux, without requiring an explicit spectral expansion or a separate pointwise reconstruction of the normal derivative \cite{Westhuis2001}. The recovered flux supplies the ALE surface equations and satisfies a global discrete flux balance.

The resulting coupling gives a direct relation between the elliptic flux balance and the geometric rate of change of fluid volume. We prove semi-discrete volume conservation under compatible mesh motion, surface projection, and quadrature. An optional sixth-order Kreiss--Oliger-type surface dissipation \cite{KreissOliger1973} is combined with a projection onto the tangent space of the constant-volume manifold. The projected elevation correction makes no contribution to the instantaneous volume rate, allowing the stabilization to retain the same semi-discrete balance. A fourth-order Runge--Kutta method, with the mesh updated at every stage, explicitly advances the surface variables. The volume result concerns the semi-discrete equations, while exact volume preservation at finite time steps is not generally guaranteed on a moving domain.

Our method combines several advantages within a unified ALE continuous finite element framework. It treats the fully nonlinear free-surface equations directly on the evolving physical domain, accommodates non-flat bottoms and prescribed boundary motion, and retains a high-order spatial discretization together with the discrete volume-balance structure described above. High-order spectral-element discretizations provide a complementary route to accurate fully nonlinear wave simulation on geometrically flexible meshes \cite{EngsigKarupEtAl2016}; recent work has also demonstrated a scalable three-dimensional spectral/$hp$-element implementation \cite{EskilssonAbediEngsigKarup2025}. The use of an interior finite element formulation also provides geometric flexibility for basin-type configurations, in the spirit of the field-based approach adopted by Westhuis \cite{Westhuis2001} for nonlinear waves in a hydrodynamic model test basin. The present method further couples this physical-domain discretization with a residual-based Dirichlet--Neumann flux and stage-wise ALE evolution. Compared with the DG method of Van der Vegt and Tomar \cite{van2005discontinuous}, which was developed and analyzed for linear free-surface waves, the present formulation addresses fully nonlinear surface evolution on a moving geometry. Compared with the conservative high-order space--time DG method of Van der Vegt and Xu \cite{vanDerVegtXu2007}, it avoids a global space--time discretization and the associated nonlinear algebraic solver: at each Runge--Kutta stage, only the interior elliptic problem is solved on the instantaneous geometry, followed by acompatin explicit update of the surface variables. This yields a comparatively simple and modular high-order framework while retaining geometric flexibility and the conservation structure of the semi-discrete formulation.

The remainder of this paper is organized as follows. Section~\ref{sec:governing} introduces the governing equations and the free-surface formulation. Section~\ref{sec:numerical_method} presents the Laplace solver, the volume-conserving ALE surface method with optional dissipation, and the time discretization. Section~\ref{sec:numerical_results} examines spatial accuracy and nonlinear waves in fixed- and moving-boundary configurations. Section~\ref{sec:conclusion} summarizes the findings and future directions.

\section{Governing equations and free-surface formulation}\label{sec:governing}
\subsection{Laplace equation and boundary conditions}
\mbox{}\par
We consider two-dimensional inviscid, incompressible and irrotational
free-surface flows. Let $\Omega(t)$ denote the time-dependent fluid domain,
whose upper boundary is the free surface
\[
    \Gamma_S(t)=\{(x,y):y=\eta(x,t)\}.
\]
The fluid velocity is expressed in terms of the velocity potential $\phi$ as
\[
    \boldsymbol{u}=\nabla\phi.
\]
The incompressibility condition then leads to the Laplace equation
\begin{equation}
    \Delta\phi=0
    \qquad\text{in }\Omega(t).
    \label{eq:laplace}
\end{equation}

On the solid boundary $\Gamma_B(t)$, the impermeability condition is written as
\begin{equation}
    \frac{\partial\phi}{\partial n}
    =
    \boldsymbol V_B \cdot \boldsymbol n,
    \qquad \text{on }\Gamma_B(t),
    \label{eq:solid_bc}
\end{equation}
where $\boldsymbol{V}_B$ denotes the prescribed boundary velocity. In particular, $\boldsymbol V_B=\boldsymbol 0$ for a fixed impermeable
boundary.

The free surface is a material boundary. Introducing the level-set
function
\begin{equation}
    F(x,y,t)=y-\eta(x,t),
\end{equation}
the kinematic condition is expressed as
\begin{equation}
    \frac{DF}{Dt}
    =
    F_t+\boldsymbol{u}\cdot\nabla F=0.
\end{equation}
Since
\[
    F_t=-\eta_t,
    \qquad
    \nabla F=(-\eta_x,1)^T,
\]
the two-dimensional kinematic free-surface condition becomes
\begin{equation}
    \eta_t
    =
    \phi_y-\eta_x\phi_x,
    \qquad\text{on }\Gamma_S(t).
    \label{eq:kinematic}
\end{equation}

The dynamic free-surface condition follows from the Bernoulli equation.
Taking the atmospheric pressure as the reference pressure, we
obtain
\begin{equation}
    \phi_t
    +
    \frac{1}{2}
    \left(
    \phi_x^2+\phi_y^2
    \right)
    +g\eta
    =0,
    \qquad\text{on }\Gamma_S(t).
    \label{eq:dynamic}
\end{equation}

The boundary of the computational domain is decomposed into the free
surface $\Gamma_S(t)$, solid boundaries $\Gamma_B(t)$, and, when present,
periodic boundaries $\Gamma_P(t)$. Periodic boundary pairs can be treated as interior interfaces
in the finite element discretization and therefore omitted from the
subsequent discussion of physical boundary conditions.

\subsection{Free-surface formulation}

To formulate the free-surface dynamics in terms of surface variables, we follow Zakharov's formulation \cite{Zakharov1968} and use the free-surface elevation together with the trace of the velocity potential as the primary unknowns. The resulting evolution equations are expressed through the Dirichlet--Neumann operator, following the surface representation of Craig and Sulem \cite{CraigSulem1993}. To accommodate moving solid boundaries, we retain the prescribed solid-boundary flux in the interior potential problem. We then define the surface potential by
\begin{equation}
    \psi(x,t)
    =
    \phi(x,\eta(x,t),t).
    \label{eq:psi_definition}
\end{equation}

Differentiating \eqref{eq:psi_definition} with respect to $x$ and $t$
gives
\begin{align}
    \psi_x
    &=
    \phi_x+\eta_x\phi_y,
    \label{eq:psi_x}\\
    \psi_t
    &=
    \phi_t+\eta_t\phi_y,
    \label{eq:psi_t_chain}
\end{align}
where the derivatives of $\phi$ are evaluated at the free surface.

Given \(\eta\), \(\psi\), and the prescribed solid-boundary flux \(g_N:=\boldsymbol V_B\cdot\boldsymbol n\), let \(\phi\) be the solution of the corresponding Laplace problem in \(\Omega(t)\). The Dirichlet--Neumann operator is then defined by
\begin{equation}
    G(\eta;g_N)\psi
    =
    \phi_y-\eta_x\phi_x, \qquad\text{on  }  y=\eta(x,t).
    \label{eq:dno_definition}
\end{equation}

In the following, $G(\eta)\psi$ abbreviates $G(\eta;g_N)\psi$ when
the prescribed solid-boundary data are understood.
Since the outward unit normal to the free surface is
\begin{equation}
    \boldsymbol{n}
    =
    \frac{(-\eta_x,1)^T}
    {\sqrt{1+\eta_x^2}},
\end{equation}
one has
\begin{equation}
    G(\eta)\psi
    =
    \sqrt{1+\eta_x^2}
    \frac{\partial\phi}{\partial n}.
\end{equation}

Using \eqref{eq:psi_x} and $G(\eta)$, the velocity components on the free surface
can be expressed in terms of $\psi_x$ and $G(\eta)\psi$ as
\begin{equation}
    \phi_x
    =
    \frac{\psi_x-\eta_xG(\eta)\psi}{1+\eta_x^2},
    \qquad
    \phi_y
    =
    \frac{G(\eta)\psi+\eta_x\psi_x}{1+\eta_x^2}.
\end{equation}
Substituting these relations into the kinematic and dynamic conditions \eqref{eq:kinematic}--\eqref{eq:dynamic} yields the closed surface evolution system
\begin{align}
    \eta_t
    &=
    G(\eta)\psi,
    \label{eq:eta_evolution}\\
    \psi_t
    &=
    -g\eta
    -\frac{1}{2}\psi_x^2
    +\frac{1}{2}
    \frac{\left(G(\eta)\psi+\eta_x\psi_x\right)^2}
    {1+\eta_x^2}.
    \label{eq:psi_evolution}
\end{align}

Thus, for a given surface state \((\eta,\psi)\), the velocity potential \(\phi\) is obtained by solving the Laplace problem in the instantaneous fluid domain. The resulting Dirichlet--Neumann flux then determines the right-hand sides of the surface evolution equations. This naturally
leads to a coupled bulk--surface formulation: the Laplace equation is
discretized in the fluid domain, while the free-surface variables are
evolved in the associated trace space. Section~\ref{sec:numerical_method}
first constructs the interior solver and then derives the ALE surface
discretization and its volume balance.

\section{Numerical method}\label{sec:numerical_method}
Given the free-surface state and the current mesh, we first solve the
interior Laplace problem. The resulting finite element residual defines
the discrete Dirichlet--Neumann flux used in the ALE surface evolution.
We establish the semi-discrete volume balance, extend it to include
optional surface dissipation, and then describe the Runge--Kutta time
integration with a mesh update at each stage.

\subsection{Finite element solution of the Laplace problem}\label{sec:laplace_discrete}
For each time $t$, let $\mathcal{T}_h(t)$ be a conforming triangular
partition of the time-dependent fluid domain $\Omega(t)$. We introduce
the Sobolev space
\begin{equation}
    H^1(\Omega(t))
    =
    \left\{
    v\in L^2(\Omega(t)):
    \nabla v\in [L^2(\Omega(t))]^2
    \right\}.
\end{equation}
The corresponding conforming finite element space is defined by
\begin{equation}
    V_h^k(t)
    =
    \left\{
    v_h\in H^1(\Omega(t)):
    v_h|_K
    =
    \widehat{v}_h\circ\boldsymbol{F}_K^{-1},
    \quad
    \widehat{v}_h\in P^k(\widehat{K}),
    \quad
    \forall K\in\mathcal{T}_h(t)
    \right\},
    \label{eq:fem_space}
\end{equation}
where $\widehat{K}$ denotes the reference triangle and
$\boldsymbol{F}_K:\widehat{K}\rightarrow K$ is the mapping from the
reference element to the physical element.
We write $\Omega_h(t)=\bigcup_{K\in\mathcal{T}_h(t)}K$ for the domain
represented by the mesh. In discrete integrals, $\Omega(t)$ and its boundary
denote these represented geometries unless the distinction is needed.
Each element map is assumed to have a positive Jacobian, and the discrete free surface is assumed to remain a single-valued graph over the horizontal coordinate.

For higher-order approximations, we use isoparametric elements.
The geometry map $\boldsymbol{F}_K$ and the finite element approximation
use the same polynomial degree; see, e.g., \cite{BrennerScott2008}.

To impose the free-surface Dirichlet data, we introduce the spaces
\begin{align}
    V_0(t)
    &=
    \left\{
    v\in H^1(\Omega(t)):
    v|_{\Gamma_S(t)}=0
    \right\},\\
    V_\psi(t)
    &=
    \left\{
    v\in H^1(\Omega(t)):
    v|_{\Gamma_S(t)}=\psi
    \right\}.
\end{align}

Multiplying the Laplace equation by a test function
$v\in V_0(t)$ and applying integration by parts yields the following
variational formulation: Find $\phi\in V_\psi(t)$ such that
\begin{equation}
    a(\phi,v)
    =
    \ell(v),
    \qquad
    \forall v\in V_0(t),
    \label{eq:weak_laplace}
\end{equation}
where
\begin{align}
    a(\phi,v)
    &=
    \int_{\Omega(t)}
    \nabla\phi\cdot\nabla v\,d\boldsymbol{x},\\
    \ell(v)
    &=
    \int_{\Gamma_B(t)}
    g_N v\,ds.
\end{align}
Here, \(g_N=\boldsymbol V_B\cdot\boldsymbol n\) denotes the prescribed normal flux on the solid boundary and vanishes on fixed impermeable boundaries.

For a prescribed discrete surface potential $\psi_h$ belonging to the trace of
$V_h^k(t)$, we restrict the variational problem to the following
finite-dimensional spaces:
\begin{align}
    V_{h,0}^k(t)
    &=
    V_h^k(t)\cap V_0(t),\\
    V_{h,\psi}^k(t)
    &=
    \left\{
    v_h\in V_h^k(t):
    v_h|_{\Gamma_S(t)}=\psi_h
    \right\}.
\end{align}
The Galerkin approximation then reads: Find
$\phi_h\in V_{h,\psi}^k(t)$ such that
\begin{equation}
    a(\phi_h,v_h)
    =
    \ell(v_h),
    \qquad
    \forall v_h\in V_{h,0}^k(t).
    \label{eq:fem_laplace}
\end{equation}

Expanding $\phi_h$ in terms of the finite element basis functions leads
to a sparse linear system for the nodal values of the velocity potential.
The Dirichlet condition $\phi_h=\psi_h$ is imposed strongly on the
free-surface degrees of freedom. Although the Galerkin residual vanishes on test functions with zero free-surface trace, its action on functions with nonzero trace contains the surface flux information used below to construct the discrete Dirichlet--Neumann operator.

\subsection{Volume-conserving ALE free-surface method}\label{sec:ale}
We now couple the interior Laplace solution with the ALE free-surface evolution. The global balance of the residual-based Dirichlet--Neumann flux, together with a compatible ALE surface projection, provides the basis for semi-discrete volume conservation. This property also motivates the construction of the volume-preserving surface dissipation introduced below.

\subsubsection{Discrete Dirichlet--Neumann operator}\label{sec:discrete_dno}
Let $\gamma_S:H^1(\Omega(t))\rightarrow H^{1/2}(\Gamma_S(t))$
denote the trace operator on the free surface. We define the discrete
trace space induced by $V_h^k(t)$ as
\begin{equation}
    W_h^k(t)
    :=
    \gamma_S\bigl(V_h^k(t)\bigr)
    =
    \left.V_h^k(t)\right|_{\Gamma_S(t)}.
    \label{eq:trace_space}
\end{equation}
Let $\{\ell_a\}_{a=1}^{N_\Gamma}$ denote the corresponding nodal basis,
where $N_\Gamma$ is the number of free-surface degrees of freedom. The
discrete Dirichlet--Neumann flux is represented as an element of
$W_h^k(t)$,
\begin{equation}
    G_h(\eta_h;g_N)\psi_h
    =
    \sum_{a=1}^{N_\Gamma} G_a\ell_a
    \in W_h^k(t),
    \label{eq:discrete_dno_expansion}
\end{equation}
where $\{G_a\}_{a=1}^{N_\Gamma}$ are the coefficients to be determined.
As in the continuous formulation, we write $G_h(\eta_h)\psi_h$ when the prescribed solid-boundary flux $g_N$ is understood.

Let $I_S(t)$ denote the horizontal projection of the free surface:
\begin{equation}
    I_S(t)
    :=
    \pi_x\bigl(\Gamma_S(t)\bigr)
    =
    \left\{
    x\in\mathbb{R}:
    (x,\eta(x,t))\in\Gamma_S(t)
    \right\}.
    \label{eq:surface_interval}
\end{equation}
The free surface admits the graph representation
\begin{equation}
    \Gamma_S(t)
    =
    \left\{
    (x,\eta(x,t)):
    x\in I_S(t)
    \right\}.
\end{equation}

The continuous DNO includes the graph factor that converts the normal
flux integrated with $ds$ into a flux integrated with $dx$. Accordingly,
we use the horizontal measure in the discrete definition below.
For any $\mu_h\in W_h^k(t)$, choose an extension
$E_h\mu_h\in V_h^k(t)$ whose free-surface trace equals $\mu_h$. We define
the discrete Dirichlet--Neumann flux variationally by
\begin{equation}
    \int_{I_S(t)}\mu_h\,G_h(\eta_h;g_N)\psi_h\,dx
    :=
    a(\phi_h,E_h\mu_h)-\ell(E_h\mu_h).
    \label{eq:dno_residual_definition}
\end{equation}
This definition is independent of the extension: the difference between
two extensions has zero free-surface trace, and its residual vanishes by
\eqref{eq:fem_laplace}. Thus the recovered flux is characterized by its
variational action on trace functions.
For each basis function $\ell_a$, let $N_a=E_h\ell_a$. Substituting
\eqref{eq:discrete_dno_expansion} into
\eqref{eq:dno_residual_definition} then gives
\begin{equation}
    \sum_{b=1}^{N_\Gamma}
    \left(
    \int_{I_S(t)}\ell_a\ell_b\,dx
    \right)G_b
    =
    \int_{\Omega(t)}
    \nabla\phi_h\cdot\nabla N_a\,d\boldsymbol{x}
    -
    \int_{\Gamma_B(t)}
    g_N N_a\,ds .
    \label{eq:dno_projection}
\end{equation}
Introducing the free-surface mass matrix
$M_\Gamma=M_\Gamma(t)\in\mathbb{R}^{N_\Gamma\times N_\Gamma}$
with entries
\begin{equation}
    (M_\Gamma)_{ab}
    =
    \int_{I_S(t)}\ell_a\ell_b\,dx,
    \label{eq:surface_mass}
\end{equation}
we define the residual vector $\boldsymbol{r}_\Gamma$ by
\[
    (\boldsymbol{r}_\Gamma)_a
    =
    a(\phi_h,N_a)-\ell(N_a).
\]
The coefficient vector
$\boldsymbol{G}=(G_1,\ldots,G_{N_\Gamma})^T$
of the discrete Dirichlet--Neumann flux is then determined by
\begin{equation}
    M_\Gamma\boldsymbol{G}
    =
    \boldsymbol{r}_\Gamma.
    \label{eq:discrete_dno}
\end{equation}

The resulting flux is consistent with the interior variational
formulation and requires no separate pointwise reconstruction of
$\nabla\phi_h$. We next couple this discrete flux with the ALE
free-surface equations and establish the corresponding
volume-conservation property.

\subsubsection{Volume-conserving ALE formulation}\label{sec:ale_volume}
We describe the evolving free surface using a horizontal ALE map.
Let $\xi$ denote the reference coordinate and define
\begin{equation}
    x=X(\xi,t).
\end{equation}
The corresponding horizontal mesh velocity is defined by
\begin{equation}
    w(x,t)
    :=
    \left.\frac{\partial X(\xi,t)}{\partial t}\right|_{\xi}.
\end{equation}

For a sufficiently smooth scalar function $f(x,t)$, the time derivative
at fixed reference coordinate $\xi$ satisfies
\begin{equation}
    \left.\frac{d f}{dt}\right|_{\xi}
    =
    \frac{\partial f}{\partial t}
    +
    w\frac{\partial f}{\partial x}.
    \label{eq:ALE_derivative}
\end{equation}

Applying \eqref{eq:ALE_derivative} to the surface equations
\eqref{eq:eta_evolution} and \eqref{eq:psi_evolution} gives
\begin{align}
    \left.\frac{d\eta}{dt}\right|_{\xi}
    &=
    G(\eta)\psi+w\eta_x,
    \label{eq:ALE_eta}\\
    \left.\frac{d\psi}{dt}\right|_{\xi}
    &=
    -g\eta
    -\frac{1}{2}\psi_x^2
    +
    \frac{1}{2}
    \frac{
        \left(G(\eta)\psi+\eta_x\psi_x\right)^2
    }{
        1+\eta_x^2
    }
    +w\psi_x.
    \label{eq:ALE_psi}
\end{align}

We approximate the free-surface variables in the trace space
$W_h^k(t)$ defined in \eqref{eq:trace_space}:
\begin{equation}
    \eta_h(x,t)
    =
    \sum_{a=1}^{N_\Gamma}\eta_a(t)\ell_a(x,t),
    \qquad
    \psi_h(x,t)
    =
    \sum_{a=1}^{N_\Gamma}\psi_a(t)\ell_a(x,t).
    \label{eq:surface_expansion}
\end{equation}

Their horizontal derivatives are evaluated elementwise as
\begin{equation}
    \eta_{h,x}
    =
    \sum_{a=1}^{N_\Gamma}
    \eta_a\frac{\partial\ell_a}{\partial x},
    \qquad
    \psi_{h,x}
    =
    \sum_{a=1}^{N_\Gamma}
    \psi_a\frac{\partial\ell_a}{\partial x}.
    \label{eq:surface_derivatives}
\end{equation}

Let $w_h$ denote the discrete horizontal mesh velocity
induced by the ALE mesh motion.
Using the discrete Dirichlet--Neumann flux from
Section~\ref{sec:discrete_dno}, we define the
free-surface right-hand sides as
\begin{align}
    R_{\eta,h}
    &=
    G_h(\eta_h)\psi_h
    +
    w_h\eta_{h,x},
    \label{eq:R_eta}\\
    R_{\psi,h}
    &=
    -g\eta_h
    -\frac{1}{2}\psi_{h,x}^2
    +
    \frac{1}{2}
    \frac{
        \left(
            G_h(\eta_h)\psi_h
            +\eta_{h,x}\psi_{h,x}
        \right)^2
    }{
        1+\eta_{h,x}^2
    }
    +
    w_h\psi_{h,x}.
    \label{eq:R_psi}
\end{align}

The semi-discrete free-surface equations are then obtained by projecting
$R_{\eta,h}$ and $R_{\psi,h}$ onto $W_h^k(t)$. Since the trace basis functions are transported by the ALE mapping,
they satisfy
\[
    \partial_t\ell_a
    +
    w_h\partial_x\ell_a
    =
    0.
\]
Consequently, the ALE time derivatives of the discrete surface
variables are given by
\[
    \dot\eta_h
    =
    \sum_{a=1}^{N_\Gamma}\dot\eta_a\ell_a,
    \qquad
    \dot\psi_h
    =
    \sum_{a=1}^{N_\Gamma}\dot\psi_a\ell_a.
\]

Projecting the discrete right-hand sides
$R_{\eta,h}$ and $R_{\psi,h}$ onto $W_h^k(t)$
in the $L^2(I_S(t))$ sense yields
\begin{align}
    \int_{I_S(t)}
    \ell_a\dot\eta_h\,dx
    &=
    \int_{I_S(t)}
    \ell_a R_{\eta,h}\,dx,
    \label{eq:eta_projection}\\
    \int_{I_S(t)}
    \ell_a\dot\psi_h\,dx
    &=
    \int_{I_S(t)}
    \ell_a R_{\psi,h}\,dx,
    \label{eq:psi_projection}
\end{align}
for $a=1,\ldots,N_\Gamma$. Let
\[
    \boldsymbol{\eta}
    =
    (\eta_1,\ldots,\eta_{N_\Gamma})^T,
    \qquad
    \boldsymbol{\psi}
    =
    (\psi_1,\ldots,\psi_{N_\Gamma})^T
\]
denote the coefficient vectors of the discrete surface variables.
Using the free-surface mass matrix
$M_\Gamma(t)$ defined in \eqref{eq:surface_mass},
the semi-discrete system can be written as
\begin{equation}
    M_\Gamma(t)\dot{\boldsymbol{\eta}}
    =
    \boldsymbol{b}_\eta,
    \qquad
    M_\Gamma(t)\dot{\boldsymbol{\psi}}
    =
    \boldsymbol{b}_\psi,
    \label{eq:semi_discrete_surface}
\end{equation}
where
\[
    (\boldsymbol{b}_\eta)_a
    =
    \int_{I_S(t)}
    \ell_a R_{\eta,h}\,dx,
    \qquad
    (\boldsymbol{b}_\psi)_a
    =
    \int_{I_S(t)}
    \ell_a R_{\psi,h}\,dx.
\]

The residual-based flux and the ALE kinematic equation together yield
the following semi-discrete volume-conservation property.

\begin{proposition}[Semi-discrete volume conservation]
\label{prop:volume_conservation}
Let $V_h(t)=|\Omega_h(t)|$ denote the area of a
nondegenerate discrete fluid domain whose free surface
admits a graph representation. Assume that constants
belong to the finite element space and that the prescribed
solid-boundary flux satisfies
$g_N=\boldsymbol V_B\cdot\boldsymbol n$.
Furthermore, assume that the discrete geometry,
ALE mesh motion, surface projection, and quadrature
are compatible with the geometric rate of change
of the discrete fluid area. Then the semi-discrete
formulation without artificial dissipation conserves
the discrete fluid area:
\begin{equation}
    \frac{dV_h}{dt}=0.
    \label{eq:discrete_volume_conservation}
\end{equation}
\end{proposition}

\begin{proof}
By the residual-based definition of the discrete
Dirichlet--Neumann flux,
\begin{equation}
    M_\Gamma\boldsymbol G
    =
    \boldsymbol r_\Gamma,
    \label{eq:volume_dno_relation}
\end{equation}
where $\boldsymbol G$ and $\boldsymbol r_\Gamma$
are defined in Section~\ref{sec:discrete_dno}.

Since the trace basis forms a partition of unity,
$\sum_{a=1}^{N_\Gamma}\ell_a=1$, we obtain
\begin{align}
    \int_{I_S(t)}
    G_h(\eta_h)\psi_h\,dx
    &=
    \boldsymbol 1^T M_\Gamma\boldsymbol G
    \nonumber\\
    &=
    \boldsymbol 1^T\boldsymbol r_\Gamma,
    \label{eq:total_free_surface_flux}
\end{align}
where $\boldsymbol 1=(1,\ldots,1)^T$.

Although the constant function cannot be used as a
test function in \eqref{eq:fem_laplace}, it is admissible
in the residual identity
\eqref{eq:dno_residual_definition}.
Taking $\mu_h=1$ and choosing $E_h1=1$, we obtain
\begin{equation}
    \int_{I_S(t)}
    G_h(\eta_h)\psi_h\,dx
    =
    a(\phi_h,1)-\ell(1)
    =
    -\int_{\Gamma_B(t)}g_N\,ds,
    \label{eq:discrete_flux_compatibility}
\end{equation}
where $a(\phi_h,1)=0$ because $\nabla 1=0$.
Hence, the global discrete flux balance is
\begin{equation}
    \int_{I_S(t)}
    G_h(\eta_h)\psi_h\,dx
    +
    \int_{\Gamma_B(t)}g_N\,ds
    =
    0.
    \label{eq:global_flux_zero}
\end{equation}

Summing \eqref{eq:eta_projection} over all free-surface
degrees of freedom and using the partition-of-unity
property gives
\begin{equation}
    \int_{I_S(t)}
    \left(
        \dot\eta_h-w_h\eta_{h,x}
    \right)\,dx
    =
    \int_{I_S(t)}
    G_h(\eta_h)\psi_h\,dx.
    \label{eq:integrated_kinematic_balance}
\end{equation}
This identity holds in the integrated sense and does
not imply pointwise equality of the integrands.

By the assumed compatibility of the discrete geometry
and ALE mesh motion, the geometric area derivative
satisfies
\begin{equation}
    \frac{dV_h}{dt}
    =
    \int_{I_S(t)}
    \left(
        \dot\eta_h-w_h\eta_{h,x}
    \right)\,dx
    +
    \int_{\Gamma_B(t)}
    \boldsymbol V_B\cdot\boldsymbol n\,ds.
    \label{eq:discrete_volume_rate}
\end{equation}
Combining \eqref{eq:integrated_kinematic_balance}
with $g_N=\boldsymbol V_B\cdot\boldsymbol n$
and \eqref{eq:global_flux_zero}, we obtain
\begin{equation}
    \frac{dV_h}{dt}
    =
    \int_{I_S(t)}
    G_h(\eta_h)\psi_h\,dx
    +
    \int_{\Gamma_B(t)}g_N\,ds
    =
    0.
\end{equation}
Thus, the semi-discrete formulation conserves
the discrete fluid volume.
\end{proof}

\begin{remark}
\label{rem:periodic_volume}
For periodic boundary conditions, paired periodic boundaries are
identified as interior interfaces in the global finite element
formulation. Their normal flux contributions cancel pairwise, and hence
they do not contribute to the global balance
\eqref{eq:global_flux_zero}. Therefore, the same discrete
volume-conservation property holds for periodic domains.
\end{remark}

\subsubsection{Volume-conserving numerical dissipation}
\label{sec:ale_dissipation}

Under the assumptions of
Proposition~\ref{prop:volume_conservation},
the semi-discrete formulation conserves fluid volume.
To suppress high-frequency oscillations during
long-time simulations on deforming meshes,
we introduce optional sixth-order Kreiss--Oliger (KO6)
dissipation \cite{KreissOliger1973} into the free-surface
equations, leaving the interior Laplace discretization
unchanged. The dissipative corrections are constructed
to preserve the semi-discrete volume balance.

To motivate the sixth-order dissipation,
consider a generic evolution equation
\begin{equation}
    q_t=\mathcal{R}(q),
\end{equation}
where $q$ denotes either $\eta$ or $\psi$.
Adding a sixth-order dissipative term gives
\begin{equation}
    q_t
    =
    \mathcal{R}(q)
    +
    \nu_6\frac{\partial^6q}{\partial x^6},
    \qquad \nu_6>0,
    \label{eq:continuous_ko6}
\end{equation}
where $\nu_6$ is the dissipation coefficient.
For a Fourier mode $q=e^{\mathrm{i}kx}$
on a periodic interval,
\begin{equation}
    \nu_6\partial_x^6q
    =
    -\nu_6 k^6q.
\end{equation}
Thus, the sixth-order derivative preferentially
damps high-frequency modes, motivating the discrete
filter introduced below.

We define the discrete filter using the nodal vector
\begin{equation}
    \boldsymbol{q}
    =
    (q_0,q_1,\ldots,q_{N_\Gamma-1})^T,
\end{equation}
where $N_\Gamma$ is the number of free-surface
degrees of freedom. The third-difference operator
$B$ is defined by
\begin{equation}
    (B\boldsymbol{q})_j
    =
    -q_j+3q_{j+1}-3q_{j+2}+q_{j+3},
    \qquad
    j=0,\ldots,N_\Gamma-4.
    \label{eq:third_difference}
\end{equation}
The operator $B^TB$ is symmetric positive semidefinite
and acts as a sixth-order-type high-frequency filter.
In particular,
\begin{equation}
    B\boldsymbol{1}=\boldsymbol{0},
\end{equation}
so constant nodal vectors are annihilated by
the unprojected filter. This property alone,
however, does not ensure preservation of the
discrete fluid volume.

Let
\begin{equation}
    M_L
    =
    \operatorname{diag}(m_0,\ldots,m_{N_\Gamma-1})
\end{equation}
be the lumped free-surface mass matrix, with $m_i>0$.
For a characteristic horizontal nodal spacing
$h_\Gamma$, define
\begin{equation}
    \kappa(t)
    =
    \frac{\nu_6}{64h_\Gamma^5}.
    \label{eq:ko6_kappa}
\end{equation}
The KO6 correction for the surface potential is
\begin{equation}
    \boldsymbol{Q}_\psi
    =
    -\kappa M_L^{-1}B^TB\boldsymbol{\psi}.
    \label{eq:ko6_psi}
\end{equation}
It satisfies
\[
    \boldsymbol{\psi}^TM_L\boldsymbol{Q}_\psi
    =
    -\kappa\|B\boldsymbol{\psi}\|_2^2
    \leq 0,
\]
showing the dissipative character of this correction
in the $M_L$ inner product.

To preserve the volume-conservation property of
Proposition~\ref{prop:volume_conservation},
the dissipation applied to the free-surface elevation
must be modified. In general, directly applying
\eqref{eq:ko6_psi} to $\boldsymbol{\eta}$ does not
produce a volume-conserving correction.

Let
\begin{equation}
    V_h=V_h(\boldsymbol{\eta},t)
\end{equation}
denote the discrete fluid volume, with the remaining
geometric data held fixed when differentiating with
respect to $\boldsymbol{\eta}$.
We define the volume gradient by
\begin{equation}
    \boldsymbol{g}_V
    :=
    \nabla_{\boldsymbol{\eta}}V_h
    =
    \frac{\partial V_h}{\partial\boldsymbol{\eta}}.
    \label{eq:volume_gradient}
\end{equation}
A perturbation $\boldsymbol{v}$ produces the
first-order volume variation
\begin{equation}
    \delta V_h
    =
    \boldsymbol{g}_V^T\boldsymbol{v}.
\end{equation}
Hence, $\boldsymbol{v}$ is tangent to the
instantaneous constant-volume manifold if
\begin{equation}
    \boldsymbol{g}_V^T\boldsymbol{v}=0.
\end{equation}

For $\boldsymbol{g}_V\ne\boldsymbol{0}$, we define
the $M_L$-orthogonal projection onto the
volume-tangent space by
\begin{equation}
    P_V\boldsymbol{v}
    =
    \boldsymbol{v}
    -
    M_L^{-1}\boldsymbol{g}_V
    \frac{
        \boldsymbol{g}_V^T\boldsymbol{v}
    }{
        \boldsymbol{g}_V^TM_L^{-1}\boldsymbol{g}_V
    }.
    \label{eq:volume_projection}
\end{equation}
By construction,
\begin{equation}
    \boldsymbol{g}_V^TP_V\boldsymbol{v}=0
    \qquad
    \forall\,\boldsymbol{v}\in\mathbb{R}^{N_\Gamma}.
    \label{eq:volume_projection_property}
\end{equation}

The KO6 correction for the free-surface elevation
is then defined by
\begin{equation}
    \boldsymbol{Q}_\eta
    =
    -\kappa
    P_VM_L^{-1}B^TBP_V\boldsymbol{\eta}.
    \label{eq:ko6_eta}
\end{equation}
Both projections are evaluated on the instantaneous
surface geometry. The rightmost projection restricts
the input to the volume-tangent space, whereas the
leftmost projection ensures that the resulting
correction is volume-conserving. Moreover, the
double-sided projection retains the nonpositive
quadratic contribution
\[
    \boldsymbol{\eta}^TM_L\boldsymbol{Q}_\eta
    =
    -\kappa\|BP_V\boldsymbol{\eta}\|_2^2
    \leq 0.
\]
The projection acts on the dissipative correction
and does not constitute a finite projection of
the surface state onto a constant-volume manifold.

\begin{proposition}[Volume conservation of the KO6 dissipation]
\label{prop:ko6_volume}
The artificial dissipation defined by
\eqref{eq:ko6_eta} makes no direct contribution
to the discrete fluid-volume rate:
\begin{equation}
    \boldsymbol{g}_V^T\boldsymbol{Q}_\eta=0.
    \label{eq:ko6_volume_neutrality}
\end{equation}
\end{proposition}

\begin{proof}
By \eqref{eq:ko6_eta} and
\eqref{eq:volume_projection_property},
\begin{align}
    \boldsymbol{g}_V^T\boldsymbol{Q}_\eta
    =
    -\kappa
    \boldsymbol{g}_V^TP_V
    M_L^{-1}B^TBP_V\boldsymbol{\eta}
    \nonumber
    =0.
\end{align}
as required.
\end{proof}

With dissipation, the semi-discrete surface system becomes
\begin{equation}
    \dot{\boldsymbol{\eta}}
    =
    M_\Gamma^{-1}\boldsymbol{b}_\eta
    +\boldsymbol{Q}_\eta,
    \qquad
    \dot{\boldsymbol{\psi}}
    =
    M_\Gamma^{-1}\boldsymbol{b}_\psi
    +\boldsymbol{Q}_\psi.
    \label{eq:dissipative_surface_system}
\end{equation}
Since $V_h$ depends on the surface geometry but not
directly on $\boldsymbol{\psi}$, the potential correction
$\boldsymbol{Q}_\psi$ introduces no direct volume-rate
contribution. Together with
\eqref{eq:ko6_volume_neutrality}, this implies that
the dissipative semi-discrete formulation also conserves
fluid volume under the assumptions of
Proposition~\ref{prop:volume_conservation}.

\begin{remark}\label{rem:dissipation}
For non-periodic free surfaces, the difference matrix
$B$ is constructed using only complete difference
windows within the physical interval. The product
$B^TB$ then defines the endpoint contributions using
these windows, without ghost nodes. For periodic
free surfaces, the difference windows wrap around
the interval, with nodal indices interpreted
periodically.
\end{remark}

In the numerical experiments, we scale the sixth-order
dissipation coefficient using a reference horizontal
mesh size and time step:
\begin{equation}
    \nu_6
    =
    \frac{\epsilon\,\Delta X_0^6}
         {64\,\Delta t_{\mathrm{ref}}},
\end{equation}
where $\epsilon$ is a dimensionless parameter controlling
the overall strength of the numerical dissipation,
$\Delta X_0$ denotes a reference horizontal mesh size,
and $\Delta t_{\mathrm{ref}}$ is a corresponding reference
time step. The coefficient $\kappa$ is then computed
from \eqref{eq:ko6_kappa} using the characteristic
horizontal nodal spacing $h_\Gamma$. The parameter
$\epsilon$ controls the strength of the filter;
the values used in the computations are specified
in Section~\ref{sec:numerical_results}.

\subsection{Time discretization and stage-wise mesh update}
\label{sec:time_discrete}

Given a prescribed mesh-motion rule that determines
the mesh geometry and velocity at each stage,
the semi-discrete surface equations in
Sections~\ref{sec:ale_volume} and
\ref{sec:ale_dissipation} can be written as
\begin{equation}
    \frac{d\boldsymbol{U}_h}{dt}
    =
    \mathcal{F}_h(\boldsymbol{U}_h,t),
    \qquad
    \boldsymbol{U}_h
    =
    \begin{pmatrix}
        \boldsymbol{\eta}\\
        \boldsymbol{\psi}
    \end{pmatrix}.
    \label{eq:semi_discrete_compact}
\end{equation}
Each evaluation of $\mathcal{F}_h$ involves the
interior Laplace solver, the residual-based DNO
recovery, the ALE surface projection, and the
optional KO6 corrections.

We integrate \eqref{eq:semi_discrete_compact}
using the classical four-stage, fourth-order
Runge--Kutta (RK4) method \cite{Butcher2016}.
Let $\Delta t=t^{n+1}-t^n$ denote the time-step
size and let $\boldsymbol{U}_h^n$ be the
approximation at time $t^n$.
The stage derivatives are
\begin{align}
    \boldsymbol{k}_1
    &=
    \mathcal{F}_h
    \left(\boldsymbol{U}_h^n,t^n\right),\\
    \boldsymbol{k}_2
    &=
    \mathcal{F}_h
    \left(
        \boldsymbol{U}_h^n
        +\frac{\Delta t}{2}\boldsymbol{k}_1,
        t^n+\frac{\Delta t}{2}
    \right),\\
    \boldsymbol{k}_3
    &=
    \mathcal{F}_h
    \left(
        \boldsymbol{U}_h^n
        +\frac{\Delta t}{2}\boldsymbol{k}_2,
        t^n+\frac{\Delta t}{2}
    \right),\\
    \boldsymbol{k}_4
    &=
    \mathcal{F}_h
    \left(
        \boldsymbol{U}_h^n
        +\Delta t\,\boldsymbol{k}_3,
        t^n+\Delta t
    \right).
\end{align}
The solution is updated by
\begin{equation}
    \boldsymbol{U}_h^{n+1}
    =
    \boldsymbol{U}_h^n
    +
    \frac{\Delta t}{6}
    \left(
        \boldsymbol{k}_1
        +2\boldsymbol{k}_2
        +2\boldsymbol{k}_3
        +\boldsymbol{k}_4
    \right).
    \label{eq:RK4}
\end{equation}

At each RK4 stage, the mesh geometry and velocity
are determined consistently from the corresponding
stage time, free-surface state, and prescribed
solid-boundary motion. The free-surface nodes follow
the stage elevation and horizontal ALE mapping,
whereas the remaining nodes are updated according
to the prescribed mesh-motion rule, with connectivity
unchanged. On the resulting stage geometry, the
interior Laplace problem is reassembled and solved,
and the discrete Dirichlet--Neumann flux and
free-surface right-hand sides are recomputed.
When KO6 dissipation is enabled, its coefficients
and volume-tangent projection are also evaluated
on the same stage geometry.

The mesh-motion rules used in the numerical
experiments are specified in
Section~\ref{sec:numerical_results}.
Each stage evaluation follows
\begin{equation}
    (\eta_h^{(s)},\psi_h^{(s)})
    \longrightarrow
    \mathcal{T}_h^{(s)}
    \longrightarrow
    \phi_h^{(s)}
    \longrightarrow
    G_h(\eta_h^{(s)})\psi_h^{(s)}
    \longrightarrow
    \mathcal{F}_h^{(s)},
    \qquad s=1,\ldots,4.
    \label{eq:RK_stage_procedure}
\end{equation}

Evaluating the spatial discretization on each
stage geometry ensures that the RK4 stages use
the corresponding semi-discrete right-hand sides.
The volume-conservation properties established
above apply to the semi-discrete formulation.
However, the discrete fluid volume generally
depends nonlinearly on the moving geometry and
is not necessarily a linear invariant of the
time-discrete state. Consequently, classical RK4
does not, in general, preserve the discrete
volume exactly at finite time steps.
The numerical experiments below assess the
resulting volume error together with the
accuracy of the computed free-surface motion.

\section{Numerical examples}\label{sec:numerical_results}
The linear-wave tests assess accuracy against an analytic solution
and compare wavemaker-generated profiles under mesh refinement.
The nonlinear tests examine free-surface evolution, mesh deformation,
and volume balance in the geometries considered below.
Most examples use rectangular tanks or related geometries. We denote
the numbers of horizontal and vertical subdivisions by $N_x$ and
$N_y$, respectively. Each cell of an $N_x\times N_y$ quadrilateral
grid is split along its upper-left to lower-right diagonal to form
two triangles.
\subsection{Common setup}
The following mesh and dissipation settings are used unless stated
otherwise.
To enhance the spatial resolution near the free surface, a nonuniform distribution is adopted in the vertical direction. More precisely, the reference grid coordinates are defined by
\begin{equation}\label{SiRj}
	s_i=\frac{i}{N_x},
    \qquad
    r_j=1-\left(1-\frac{j}{N_y}\right)^{p_v},
    \qquad
    i=0,\ldots,N_x,\quad j=0,\ldots,N_y,
\end{equation}

where $p_v$ denotes the vertical stretching parameter. For nonlinear
problems with fixed sidewalls, the corresponding physical coordinates
are given by
\begin{equation}
    x_i=x_L+s_i(x_R-x_L),
    \qquad
    y_{i,j}(t)
    =
    d(x_i)
    +
    r_j\left[\eta(x_i,t)-d(x_i)\right],
\end{equation}
where $d(x)$ describes the bottom boundary. The levels $r_j=0$ and
$r_j=1$ correspond to the bottom and free surface, respectively.
For $p_v>1$, the grid points are clustered toward the free surface;
we use $p_v=1.25$ unless stated otherwise. The linear tests use fixed
meshes, whereas the rotating-wavemaker test uses
\eqref{eq:rotating_mesh_map}.

Unless otherwise stated, the KO6 dissipation is employed with
\begin{equation}
    \epsilon = 0.05,
    \qquad
    \Delta X_0=\frac{L}{N_x},
    \qquad
    h_\Gamma=\frac{\Delta X_0}{k},
    \qquad
    \Delta t_{\mathrm{ref}}
    =
    \frac{0.02}{80}
    \left(
        \frac{\Delta X_0}{L/40}
    \right)^2,
\end{equation}
for the $P^k$ approximation, where $L$ is the horizontal length of
the reference domain. The reference step $\Delta t_{\mathrm{ref}}$
sets the scale of the dissipation coefficients, which are given by
\begin{equation}
    \nu_6
    =
    \frac{\epsilon\Delta X_0^6}
         {64\Delta t_{\mathrm{ref}}},
    \qquad
    \kappa
    =
    \frac{\nu_6}{64h_\Gamma^5}.
\end{equation}
\subsection{Linear free-surface waves}
\subsubsection{Accuracy test}\label{AT}
We consider two-dimensional linear free-surface waves governed by
\begin{align}\label{Poisson_AT}
-\Delta \phi=0,\quad \text{in }\Omega=[0,2]\times [-1,0].
\end{align}
At the flat bottom, we impose the homogeneous Neumann condition
\begin{align}\label{GammaN_AT}
\vec{n}\cdot \nabla \phi=0,\quad \text{at }y=-1.
\end{align}
The linearized boundary conditions on $\Gamma_S$ are
\begin{align}\label{GammaS_AT}
\dfrac{\PA{\phi}}{\PA{t}}+\eta=0,\quad \text{and}\quad \dfrac{\PA{\eta}}{\PA{t}}-\dfrac{\PA{\phi}}{\PA{y}}=0.
\end{align}
Periodic boundary conditions are imposed at $x=0$ and $x=2$.

The analytic solution of this problem is given by
\begin{align}\label{Phi_AT}
\phi(x,y,t)=A\omega \frac{\cosh(\kappa (y+1))}{\kappa\sinh (\kappa) }\sin (\kappa x-\omega t).
\end{align}
The corresponding elevation is
$\eta(x,t)=A\cos(\kappa x-\omega t)$. In this test, $A$ is the wave
amplitude, $\kappa$ denotes the wavenumber, and $\omega$ is the angular
frequency. With $h=1$ and $g=1$, the linear surface-gravity-wave dispersion relation is
\begin{align}\label{DR_AT}
\omega^2=\kappa \tanh(\kappa h).
\end{align}
The initial conditions for the free-surface elevation and surface potential $\psi=\phi|_{y=0}$ are
\begin{align}\label{IC_AT}
\eta(x,0)=A \cos (\kappa x),\qquad
\psi(x,0)=\frac{A\omega}{\kappa\tanh(\kappa)}\sin(\kappa x).
\end{align}

We set $A=0.02$ and $\kappa=\pi$. Both dissipation settings use
$\Delta t=\Delta t_{\mathrm{ref}}=0.025\Delta X_0^2$ and the same mesh
sequence. Tables~\ref{ta:AT} and \ref{ta:AT_no_ko6} report the errors
at $t=2.5$ for $\epsilon=0.05$ and $\epsilon=0$, respectively.
The $L^2$ and $L^\infty$ errors in $\eta$ are measured on $\Gamma_S$,
and those in $\phi$ are measured in $\Omega$.

\begin{table}[htbp]
\centering
\footnotesize
\renewcommand{\arraystretch}{1.15}
\setlength{\tabcolsep}{3pt}
\begin{tabular}{ccc*{8}{c}}
\hline
\multicolumn{2}{c}{}
& \multirow{2}{*}{$N_x\times N_y$}
& \multicolumn{4}{c}{Free-surface elevation $\eta$ on $\Gamma_S$}
& \multicolumn{4}{c}{Potential $\phi$ in $\Omega$} \\
\cline{4-11}
\multicolumn{2}{c}{}
&
& $L^2$ error & Order
& $L^\infty$ error & Order
& $L^2$ error & Order
& $L^\infty$ error & Order \\
\hline

\multicolumn{2}{c}{\multirow{5}{*}{$P^1$}}
& $20 \times 5$   & 2.815e-03 & --    & 2.835e-03 & --    & 5.275e-04 & --    & 1.451e-03 & --    \\
\multicolumn{2}{c}{}
& $40 \times 10$  & 7.252e-04 & 1.957 & 7.339e-04 & 1.950 & 1.307e-04 & 2.013 & 3.681e-04 & 1.979 \\
\multicolumn{2}{c}{}
& $80 \times 20$  & 1.827e-04 & 1.989 & 1.849e-04 & 1.989 & 3.259e-05 & 2.003 & 9.246e-05 & 1.993 \\
\multicolumn{2}{c}{}
& $160 \times 40$ & 4.575e-05 & 1.997 & 4.629e-05 & 1.998 & 8.144e-06 & 2.001 & 2.313e-05 & 1.999 \\
\multicolumn{2}{c}{}
& $320 \times 80$ & 1.144e-05 & 1.999 & 1.158e-05 & 1.999 & 2.036e-06 & 2.000 & 5.784e-06 & 1.999 \\
\hline

\multicolumn{2}{c}{\multirow{5}{*}{$P^2$}}
& $20 \times 5$   & 2.877e-05 & --    & 3.589e-05 & --    & 9.129e-06 & --    & 3.031e-05 & --    \\
\multicolumn{2}{c}{}
& $40 \times 10$  & 1.955e-06 & 3.879 & 2.625e-06 & 3.773 & 1.061e-06 & 3.105 & 3.999e-06 & 2.922 \\
\multicolumn{2}{c}{}
& $80 \times 20$  & 1.340e-07 & 3.867 & 2.021e-07 & 3.699 & 1.295e-07 & 3.035 & 5.008e-07 & 2.997 \\
\multicolumn{2}{c}{}
& $160 \times 40$ & 1.038e-08 & 3.691 & 1.740e-08 & 3.538 & 1.607e-08 & 3.010 & 6.214e-08 & 3.011 \\
\multicolumn{2}{c}{}
& $320 \times 80$ & 9.950e-10 & 3.383 & 1.688e-09 & 3.366 & 2.006e-09 & 3.003 & 7.721e-09 & 3.009 \\
\hline

\multicolumn{2}{c}{\multirow{5}{*}{$P^3$}}
& $20 \times 5$   & 1.3005e-06 & --    & 1.966e-06 & --    & 3.993e-07 & --    & 1.480e-06 & --    \\
\multicolumn{2}{c}{}
& $40 \times 10$  & 8.060e-08 & 4.012 & 1.043e-07 & 4.237 & 2.536e-08 & 3.977 & 1.025e-07 & 3.851 \\
\multicolumn{2}{c}{}
& $80 \times 20$  & 4.962e-09 & 4.022 & 5.770e-09 & 4.176 & 1.585e-09 & 4.000 & 6.714e-09 & 3.933 \\
\multicolumn{2}{c}{}
& $160 \times 40$ & 3.062e-10 & 4.018 & 3.334e-10 & 4.113 & 9.875e-11 & 4.005 & 4.283e-10 & 3.971 \\
\multicolumn{2}{c}{}
& $320 \times 80$ & 2.096e-11 & 3.869 & 2.197e-11 & 3.923 & 6.451e-12 & 3.936 & 2.828e-11 & 3.921 \\
\hline

\end{tabular}
\caption{Errors and convergence orders of the potential and free-surface elevation for linear free-surface waves at $t=2.5$ with KO6 dissipation ($\epsilon=0.05$).}
\label{ta:AT}
\end{table}

\begin{table}[htbp]
\centering
\footnotesize
\renewcommand{\arraystretch}{1.15}
\setlength{\tabcolsep}{3pt}
\begin{tabular}{ccc*{8}{c}}
\hline
\multicolumn{2}{c}{}
& \multirow{2}{*}{$N_x\times N_y$}
& \multicolumn{4}{c}{Free-surface elevation $\eta$ on $\Gamma_S$}
& \multicolumn{4}{c}{Potential $\phi$ in $\Omega$} \\
\cline{4-11}
\multicolumn{2}{c}{}
&
& $L^2$ error & Order
& $L^\infty$ error & Order
& $L^2$ error & Order
& $L^\infty$ error & Order \\
\hline

\multicolumn{2}{c}{\multirow{5}{*}{$P^1$}}
& $20 \times 5$ & 2.815e-03 & -- & 2.836e-03 & -- & 5.276e-04 & -- & 1.450e-03 & -- \\
\multicolumn{2}{c}{}
& $40 \times 10$ & 7.252e-04 & 1.957 & 7.339e-04 & 1.950 & 1.307e-04 & 2.013 & 3.680e-04 & 1.979 \\
\multicolumn{2}{c}{}
& $80 \times 20$ & 1.827e-04 & 1.989 & 1.849e-04 & 1.989 & 3.260e-05 & 2.003 & 9.246e-05 & 1.993 \\
\multicolumn{2}{c}{}
& $160 \times 40$ & 4.575e-05 & 1.997 & 4.630e-05 & 1.998 & 8.144e-06 & 2.001 & 2.313e-05 & 1.999 \\
\multicolumn{2}{c}{}
& $320 \times 80$ & 1.144e-05 & 1.999 & 1.158e-05 & 1.999 & 2.036e-06 & 2.000 & 5.784e-06 & 1.999 \\
\hline

\multicolumn{2}{c}{\multirow{5}{*}{$P^2$}}
& $20 \times 5$ & 4.768e-05 & -- & 1.051e-04 & -- & 8.569e-06 & -- & 2.816e-05 & -- \\
\multicolumn{2}{c}{}
& $40 \times 10$ & 1.209e-05 & 1.980 & 2.738e-05 & 1.940 & 1.039e-06 & 3.044 & 3.802e-06 & 2.889 \\
\multicolumn{2}{c}{}
& $80 \times 20$ & 1.268e-06 & 3.252 & 2.872e-06 & 3.253 & 1.286e-07 & 3.015 & 4.746e-07 & 3.002 \\
\multicolumn{2}{c}{}
& $160 \times 40$ & 2.481e-07 & 2.354 & 5.571e-07 & 2.366 & 1.604e-08 & 3.002 & 5.998e-08 & 2.984 \\
\multicolumn{2}{c}{}
& $320 \times 80$ & 5.527e-08 & 2.166 & 1.234e-07 & 2.175 & 2.006e-09 & 3.000 & 8.652e-09 & 2.793 \\
\hline

\multicolumn{2}{c}{\multirow{5}{*}{$P^3$}}
& $20 \times 5$ & 2.983e-06 & -- & 4.188e-06 & -- & 3.192e-07 & -- & 9.839e-07 & -- \\
\multicolumn{2}{c}{}
& $40 \times 10$ & 3.191e-07 & 3.225 & 4.408e-07 & 3.248 & 2.033e-08 & 3.972 & 7.241e-08 & 3.764 \\
\multicolumn{2}{c}{}
& $80 \times 20$ & 2.258e-08 & 3.821 & 2.944e-08 & 3.904 & 1.277e-09 & 3.993 & 4.868e-09 & 3.895 \\
\multicolumn{2}{c}{}
& $160 \times 40$ & 1.681e-09 & 3.748 & 2.283e-09 & 3.689 & 8.016e-11 & 3.994 & 3.108e-10 & 3.969 \\
\multicolumn{2}{c}{}
& $320 \times 80$ & 2.090e-10 & 3.008 & 3.579e-10 & 2.673 & 5.037e-12 & 3.992 & 2.093e-11 & 3.892 \\
\hline

\end{tabular}
\caption{Errors and convergence orders of the potential and free-surface elevation for linear free-surface waves at $t=2.5$ without KO6 dissipation ($\epsilon=0$).}
\label{ta:AT_no_ko6}
\end{table}

With KO6 dissipation, Table~\ref{ta:AT} gives approximately $(k+1)$-st
order convergence for the potential with $P^k$, $k=1,2,3$. The
elevation is approximately second order for $P^1$ and fourth order for
$P^3$; the $P^2$ rates exceed three but decrease under refinement.
Without dissipation, Table~\ref{ta:AT_no_ko6} gives the same potential
orders in $L^2$, while the elevation rates are less uniform for $P^2$
and $P^3$. KO6 has little effect on the $P^1$ errors and reduces the
elevation errors for the higher-order approximations.

\subsubsection{Linear waves generated by a wavemaker}\label{WM}

We next consider a wavemaker in the fixed domain
$\Omega=[0,10]\times[-1,0]$, using the Laplace equation and linearized
surface conditions of the previous test. Homogeneous Neumann
conditions are imposed at the bottom $y=-1$ and right wall $x=10$.
Initially, $\eta(x,0)=\psi(x,0)=0$. The interior potential is
obtained from the Laplace problem with the prescribed wall flux.
The left wall has the time-periodic condition
\begin{align}\label{GammaL_WM}
\dfrac{\partial \phi}{\partial x}=B(y+1)\cos (\omega t),\quad \text{at }x=0.
\end{align}
Since the outward unit normal at $x=0$ is $\vec n=(-1,0)$, the corresponding outward normal flux is $\vec n\cdot\nabla\phi=-B(y+1)\cos(\omega t)$.

We set $B=0.05$ and $\omega=2$ in \eqref{GammaL_WM}.
The linearization keeps the wavemaker and free-surface boundaries
at their equilibrium positions, so the mesh remains fixed.
The prescribed wall flux drives the waves through the boundary
condition. All profiles in this example include KO6 dissipation.
We first compare meshes at each polynomial degree, then compare
the polynomial degrees on each mesh.
\begin{figure}[!htbp]
\centering
\subfigure[Wave profile for $T=20$]{
\includegraphics[width=0.485\textwidth]{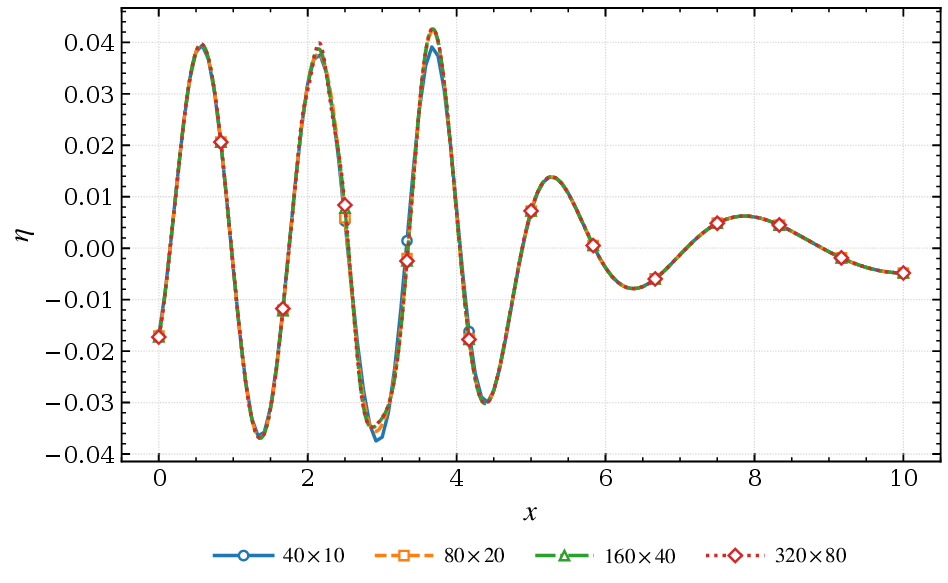}}
\subfigure[Wave profile for $T=40$]{
\includegraphics[width=0.485\textwidth]{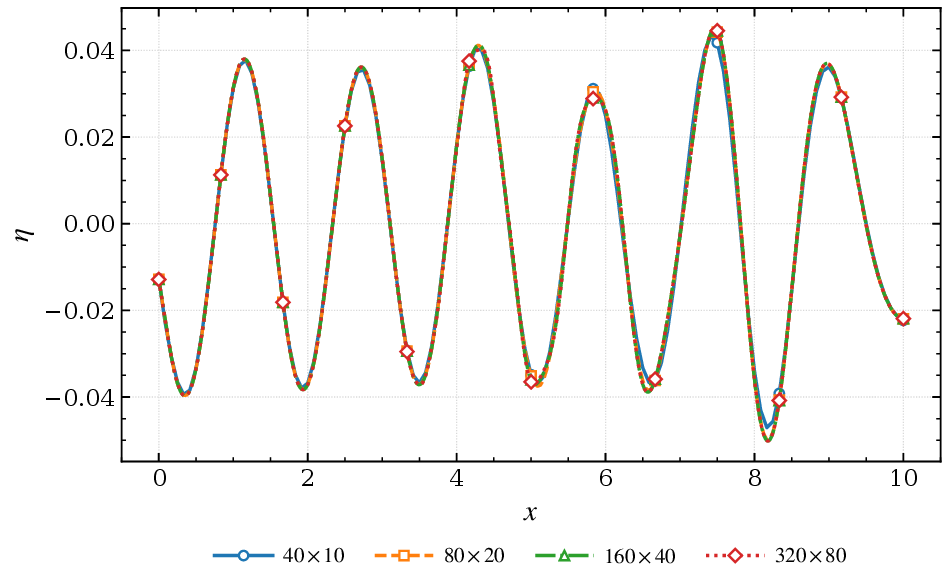}}
\subfigure[Wave profile for $T=60$]{
\includegraphics[width=0.485\textwidth]{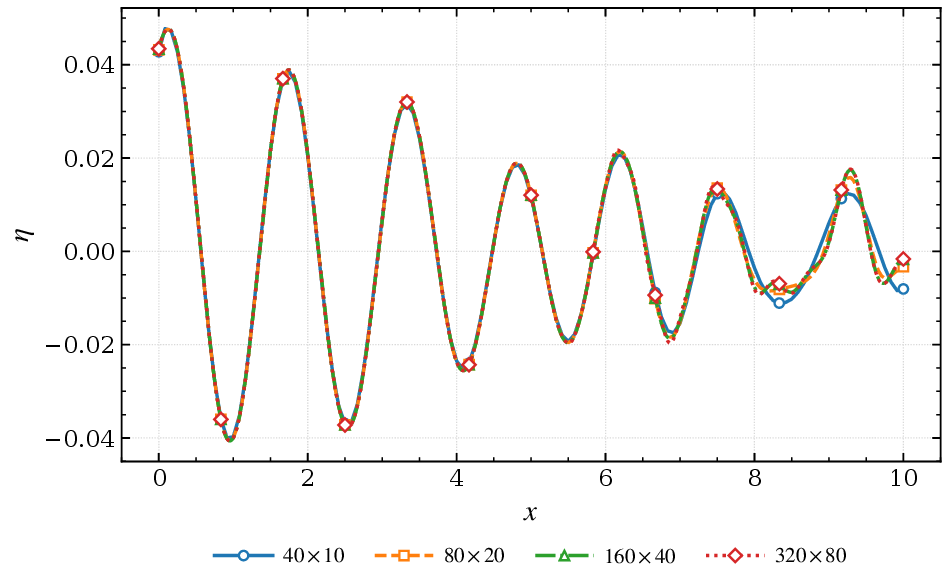}}
\subfigure[Wave profile for $T=80$]{
\includegraphics[width=0.485\textwidth]{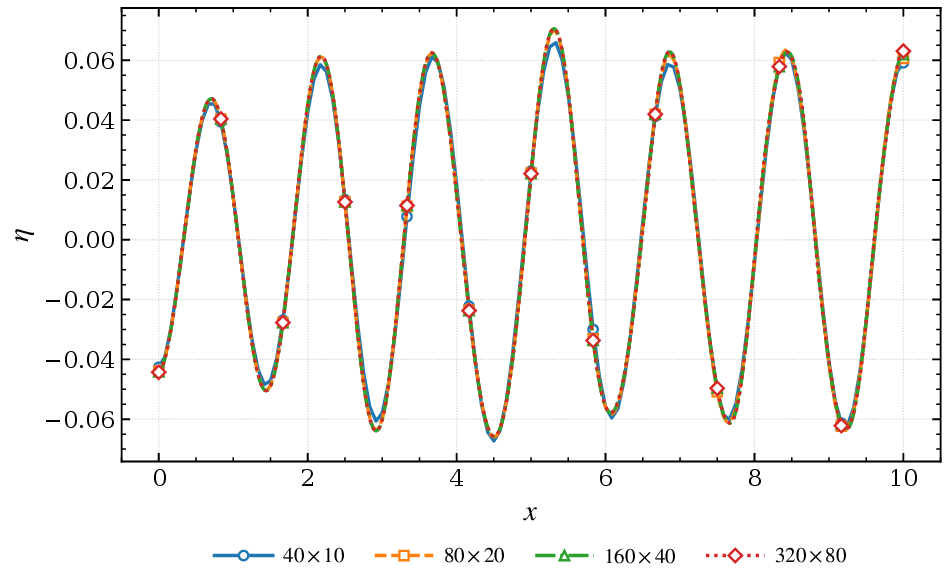}}
\caption{Wavemaker-generated free-surface profiles on different meshes using the $P^3$ approximation.}
\label{fg:WM-DN-3}
\end{figure}

\begin{figure}[!htbp]
\centering
\subfigure[Wave profile for $T=20$]{
\includegraphics[width=0.485\textwidth]{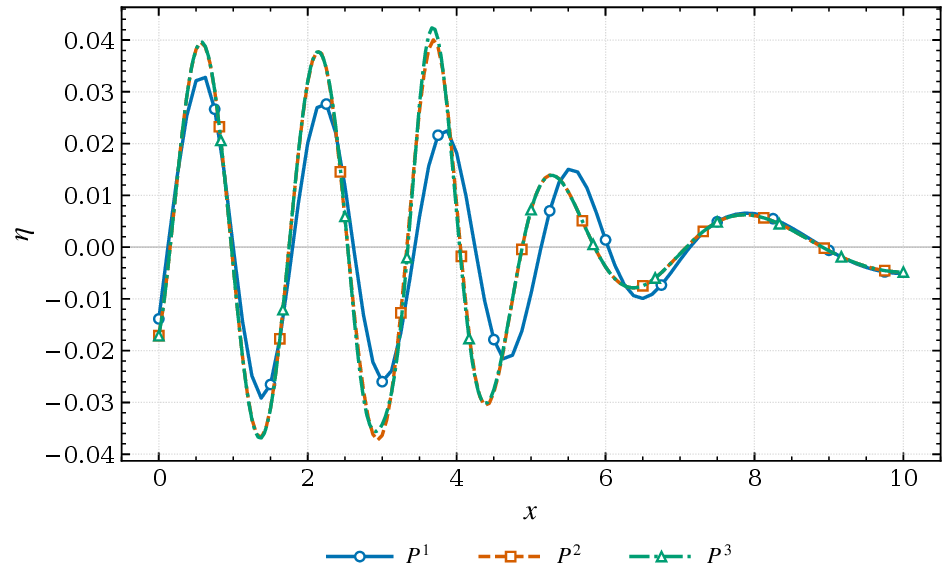}}
\subfigure[Wave profile for $T=40$]{
\includegraphics[width=0.485\textwidth]{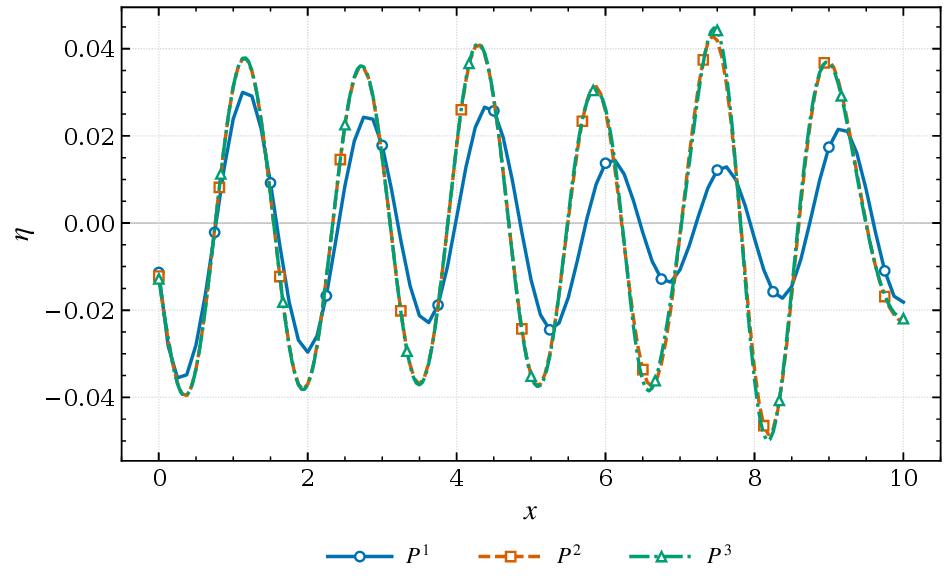}}
\subfigure[Wave profile for $T=60$]{
\includegraphics[width=0.485\textwidth]{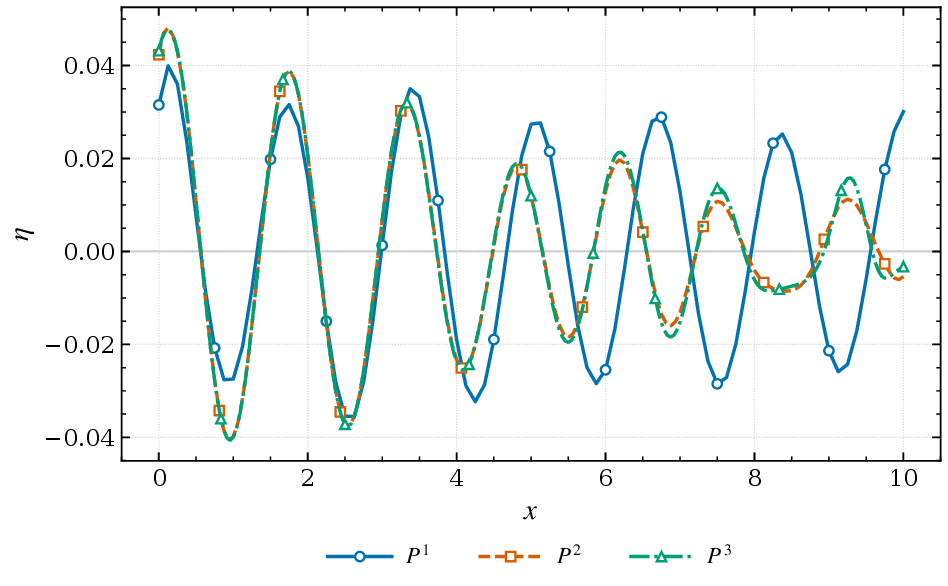}}
\subfigure[Wave profile for $T=80$]{
\includegraphics[width=0.485\textwidth]{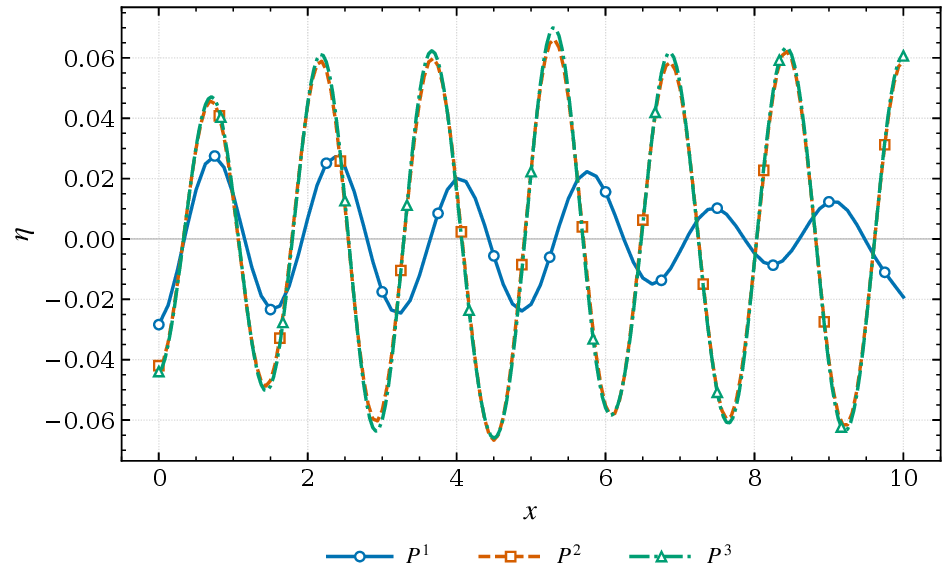}}
\caption{Wavemaker-generated free-surface profiles for the $P^1$, $P^2$, and $P^3$ approximations on an $80\times 20$ mesh.}
\label{fg:WM-DP-2}
\end{figure}
\begin{figure}[!htbp]
\centering
\subfigure[Wave profile for $T=20$]{
\includegraphics[width=0.485\textwidth]{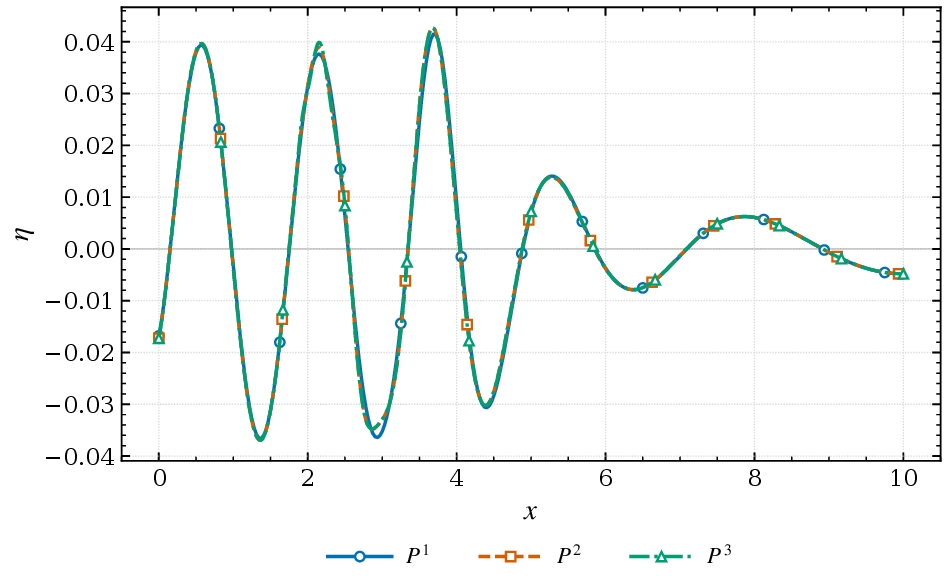}}
\subfigure[Wave profile for $T=40$]{
\includegraphics[width=0.485\textwidth]{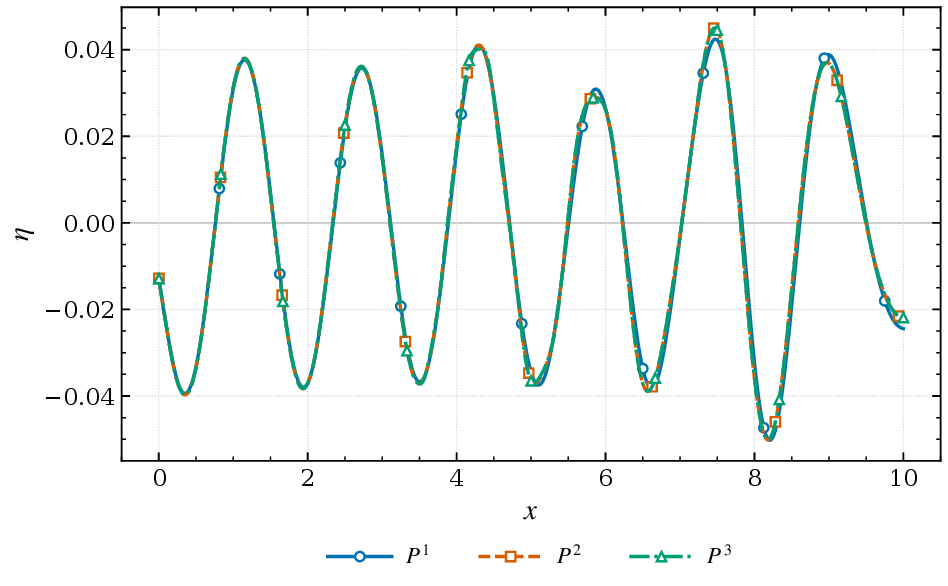}}
\subfigure[Wave profile for $T=60$]{
\includegraphics[width=0.485\textwidth]{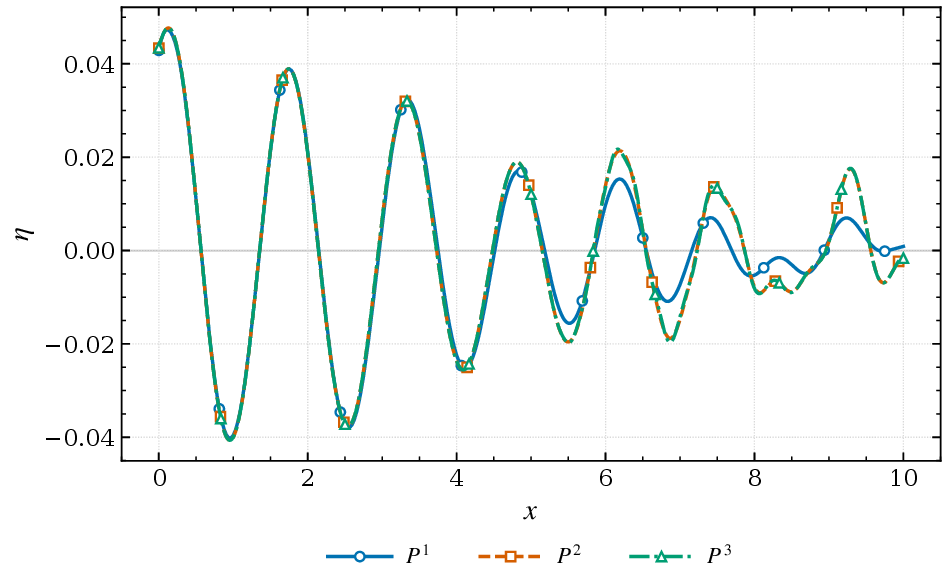}}
\subfigure[Wave profile for $T=80$]{
\includegraphics[width=0.485\textwidth]{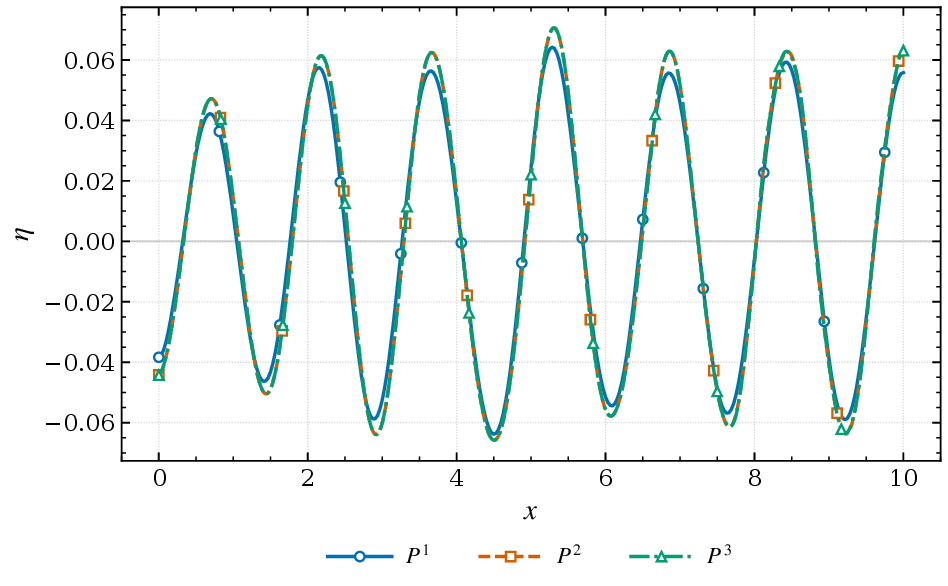}}
\caption{Wavemaker-generated free-surface profiles for the $P^1$, $P^2$, and $P^3$ approximations on a $320\times 80$ mesh.}
\label{fg:WM-DP-4}
\end{figure}

Figure~\ref{fg:WM-DN-3} shows the mesh-refinement behavior of the
$P^3$ approximation. Figures~\ref{fg:WM-DP-2} and~\ref{fg:WM-DP-4}
compare $P^1$, $P^2$, and $P^3$ on the $80\times20$ and $320\times80$
meshes, respectively. The $P^3$ profile on the $80\times20$ mesh is
already close to the corresponding profile on the $320\times80$ mesh,
whereas the $P^1$ profile changes more visibly and requires the finer
mesh to approach the apparently mesh-converged shape. This comparison
indicates that the higher-order approximation achieves comparable
resolution with substantially fewer degrees of freedom.

\subsection{Nonlinear free-surface waves}
\subsubsection{Slowly modulated wave packet}
We next consider a nonlinear wave packet governed by the Zakharov
formulation \eqref{eq:laplace}, \eqref{eq:eta_evolution}, and
\eqref{eq:psi_evolution}. 

The computational domain is given by
\[
    \Omega(t)
    =
    \{(x,y): 0\leq x\leq 2\pi,\;-1\leq y\leq \eta(x,t)\},
\]
with the initial conditions
\begin{equation}
    \eta(x,0)
    =
    0.01\exp\left(-\frac{3}{4}(x-\pi)^2\right)\cos(10x),
    \qquad
    \psi(x,0)=0.
\end{equation}
Periodic boundary conditions are imposed in the horizontal direction.
The moving upper boundary is the free surface, and the fixed bottom
satisfies the homogeneous Neumann condition $\partial_n\phi=0$.
All results in this example include KO6 dissipation.

\begin{figure}[htbp]
\centering
\begin{minipage}[t]{0.485\textwidth}
\centering
\includegraphics[width=\linewidth]{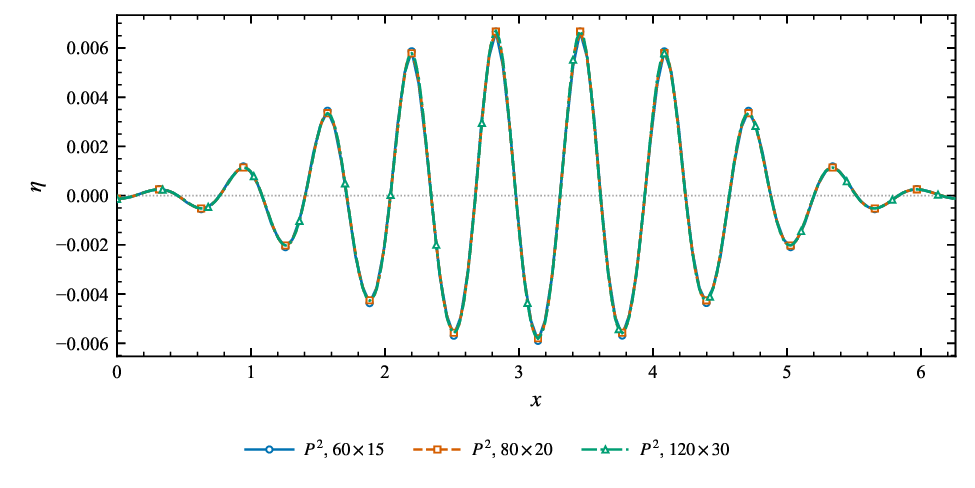}
\par\vspace{0.15em}
{\footnotesize (a) $t=5$, free surface}
\end{minipage}
\hfill
\begin{minipage}[t]{0.485\textwidth}
\centering
\includegraphics[width=\linewidth]{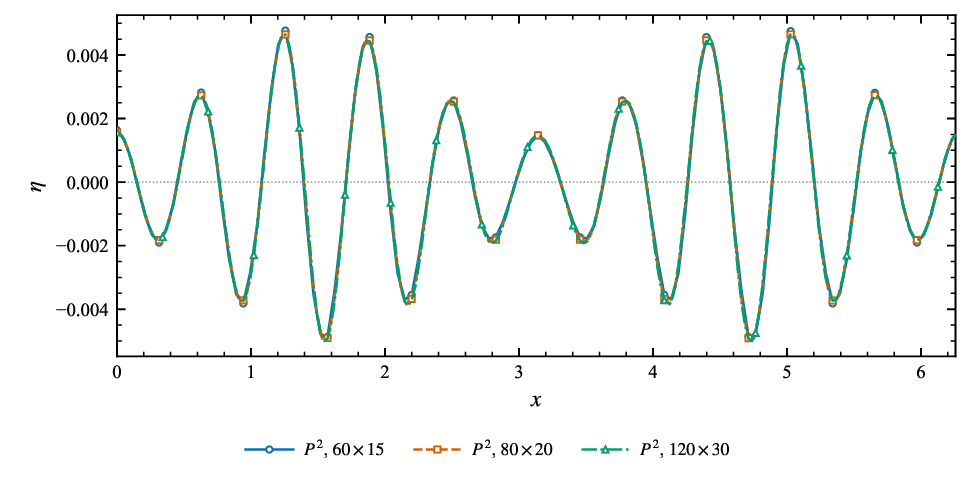}
\par\vspace{0.15em}
{\footnotesize (b) $t=10$, free surface}
\end{minipage}
\par\vspace{0.5em}
\begin{minipage}[t]{0.485\textwidth}
\centering
\includegraphics[width=\linewidth]{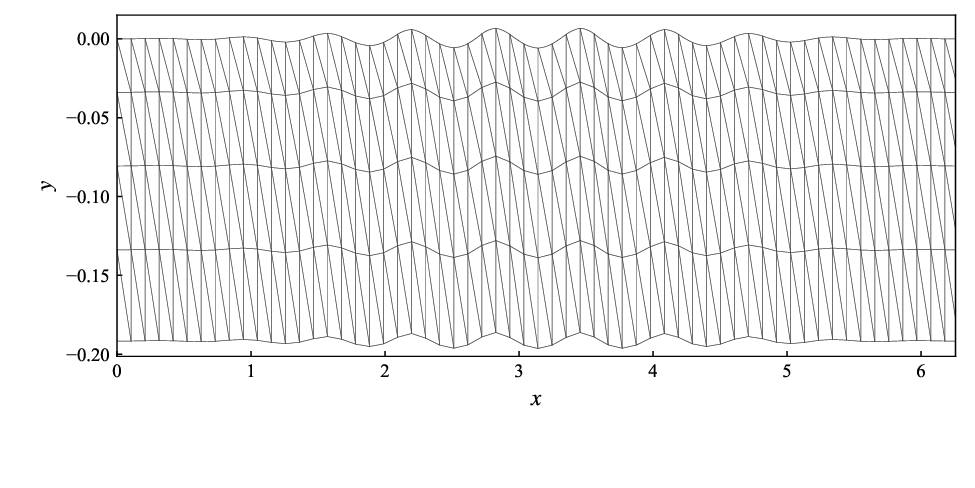}
\par\vspace{0.15em}
{\footnotesize (c) $t=5$, near-surface mesh}
\end{minipage}
\hfill
\begin{minipage}[t]{0.485\textwidth}
\centering
\includegraphics[width=\linewidth]{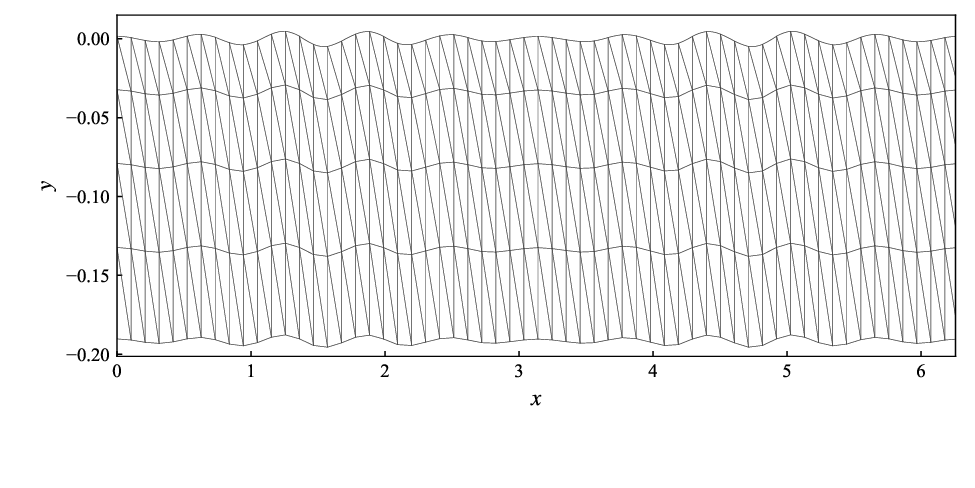}
\par\vspace{0.15em}
{\footnotesize (d) $t=10$, near-surface mesh}
\end{minipage}
\caption{Free-surface profiles and near-surface ALE meshes for the wave packet at $t=5$ (left) and $t=10$ (right). The rows show $P^2$ profiles on three meshes and the upper four mesh layers. Mesh views use the $60\times15$ grid.}
\label{fg:Pa-P2convergence}
\end{figure}

Figure~\ref{fg:Pa-P2convergence} shows that the $P^2$ profiles become
closer under mesh refinement at both $t=5$ and $t=10$. The near-surface
mesh panels show that the upper layers follow the free surface while
preserving connectivity.

\subsubsection{Rectangular basin}
To compare the nonlinear evolution with analytical references,
we consider waves in a rectangular basin governed by
\eqref{eq:laplace}, \eqref{eq:eta_evolution}, and
\eqref{eq:psi_evolution} in the domain
\[
\Omega(t)
=
\left\{
(x,y)\in \mathbb{R}^2
:
x\in [0,b],
\;
-1 < y < \eta(x,t)
\right\}.
\]
The upper boundary $y=\eta(x,t)$ is the free surface. Homogeneous
Neumann conditions enforce impermeability on the fixed solid boundaries.

We prescribe the following single-mode initial conditions:
\begin{equation}
\eta(x,0) = a \cos(\pi x), 
\qquad \phi(x,y,0)=0,
\end{equation}
where $a$ is the initial wave amplitude.

We track the normalized free-surface elevation at the basin midpoint
$x=b/2$ to compare its time history with the analytical references:
\[
\frac{\eta(b/2,t)}{a}.
\]
We set $a=0.05$ and $b=2$, and the simulation is carried out for $t\in[0,20]$.

\begin{figure}[htbp]
\centering
\includegraphics[width=0.80\textwidth]{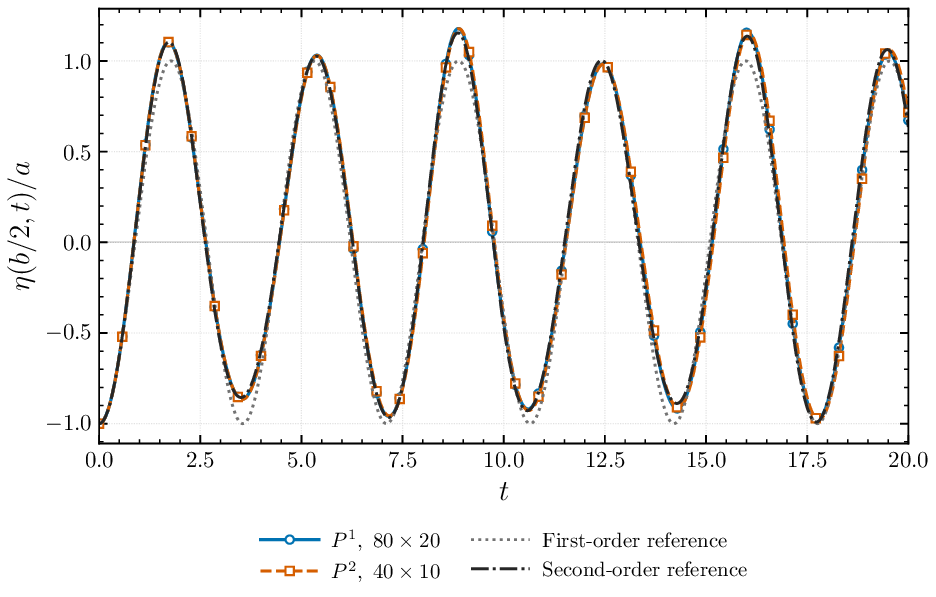}
\caption{Time history of the normalized wave elevation $\eta(b/2,t)/a$ obtained with the $P^1$ and $P^2$ approximations using KO6 dissipation, compared with the first- and second-order analytical reference solutions.}
\label{fg:Rec}
\end{figure}

Figure~\ref{fg:Rec} compares the midpoint history with the first- and
second-order analytical references of Wu and Eatock Taylor
\cite{wu1994finite}. The $P^1$ solution on the $80\times20$ mesh and
the $P^2$ solution on the $40\times10$ mesh are nearly indistinguishable
and follow the second-order reference more closely than the first-order
reference. 

\subsubsection{Basin with a bottom bump}
\label{ex:HOS_Zakharov}

We next introduce a bottom bump into the rectangular-basin problem.
The governing equations, impermeable solid boundaries, and initial
conditions are retained, with $a=0.05$ and $b=2$ over $0\leq t\leq20$. 
The computational domain is defined by
\begin{equation}
    \Omega(t)
    =
    \left\{
    (x,y)\in\mathbb{R}^2:
    x\in [0,b],\quad d(x) < y< \eta(x,t)\right \}.
\end{equation}
The bottom profile is defined by
\begin{equation}
    d(x)=
    \begin{cases}
        -2(x-\frac{1}{2})(x-\frac{3}{2})-1, & \frac{1}{2}< x< \frac{3}{2},\\
        -1, & \text{otherwise}.
    \end{cases}
\end{equation}

For this example, the mesh map uses the prescribed bottom profile
$d(x)$ and $p_v=1$ in \eqref{SiRj}. The reference vertical levels are
therefore uniform, and the physical grid follows the bump and free surface.
\begin{figure}[htbp]
\centering
\subfigure{
\includegraphics[width=0.80\textwidth]{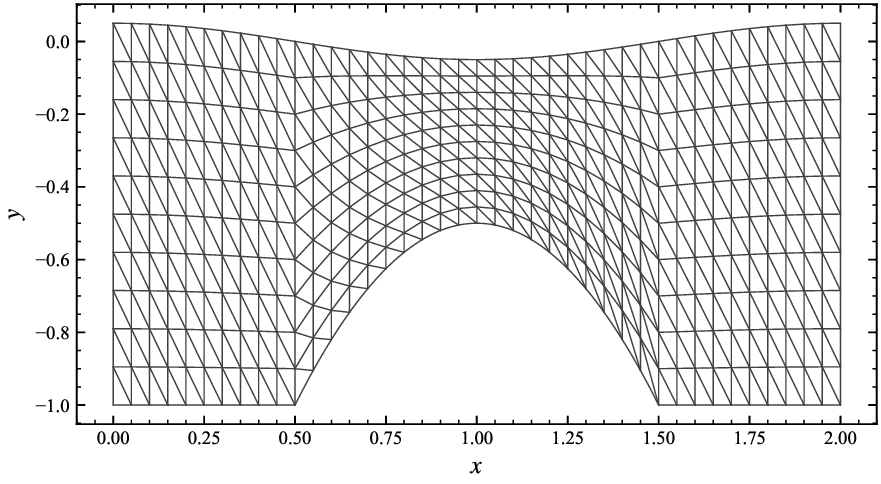}}
\subfigure{
\includegraphics[width=0.80\textwidth]{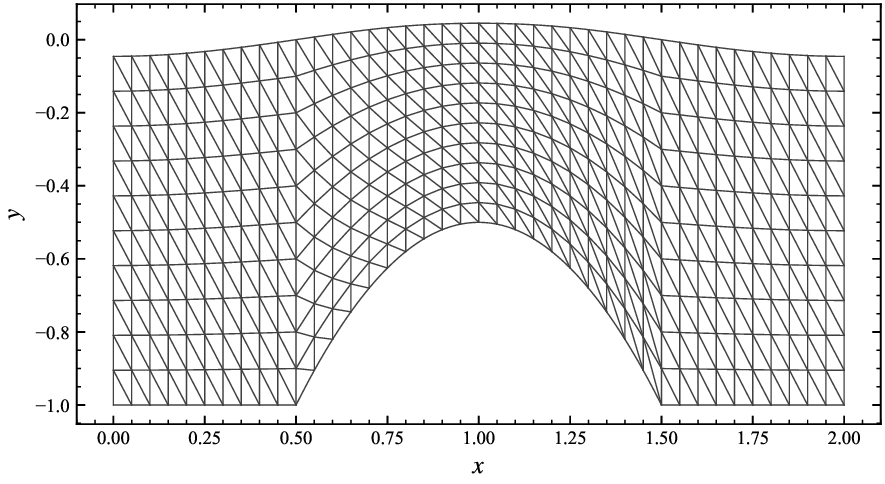}}
\caption{Computational $40\times10$ mesh for the $P^2$ approximation with KO6 dissipation in the basin with a bump at $t=0$ (top) and $t=20$ (bottom).}
\label{fg:Bump-Mesh}
\end{figure}

The mesh places nodes at $x=1/2$ and $x=3/2$, where the bottom slope
changes discontinuously. For $P^2$, quadratic isoparametric mappings
represent the adjacent curved bottom and free-surface edges. Figure~\ref{fg:Bump-Mesh}
shows the initial and deformed meshes.

\begin{figure}[htbp]
\centering
\includegraphics[width=0.80\textwidth]{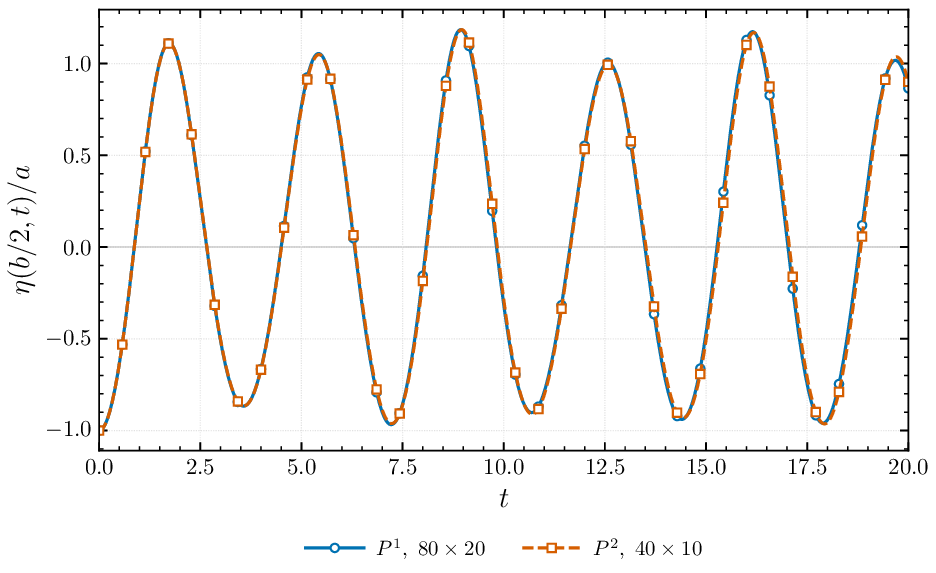}
\caption{Normalized midpoint elevation $\eta(b/2,t)/a$ in the basin with a bump for the $P^1$ and $P^2$ approximations with KO6 dissipation.}
\label{fg:Bump-Midpoint}
\end{figure}
The midpoint histories in Figure~\ref{fg:Bump-Midpoint} are closely
consistent for the $P^1$ and $P^2$ approximations over the full time
interval. Their nearly overlapping amplitudes and phases indicate that
both polynomial degrees capture the same wave evolution as the wave
passes over the bottom bump, with no visible loss of consistency in the
presence of the curved bottom.

\subsubsection{Rotating wavemaker}
Finally, we consider nonlinear waves generated by a rotating
wavemaker, extending the preceding fixed-wall tests to a moving
solid boundary. The tank has length $L$ and still-water depth $H$. The left
wall rotates about the hinge point $(0,-H)$. We prescribe its
time-dependent slope as
\begin{equation}
    q(t)=\frac{A}{H}\sin(\omega t).
\end{equation}
Its instantaneous position is given by
\begin{equation}
    x_L(y,t)=q(t)(y+H).
\end{equation}
Consequently, the fluid domain is
\begin{equation}
    \Omega(t)
    =
    \left\{
    (x,y):
    x_L(y,t)<x<L,\;
    -H<y<\eta(x,t)
    \right\}.
\end{equation}
The wall motion imposes the nonhomogeneous Neumann condition
\[
\frac{\partial\phi}{\partial\boldsymbol n}
=-\dot{\alpha}(t)\rho,
\qquad \text{on }\Gamma_L(t),
\]
where
\[
\alpha(t)=\arctan q(t),\qquad
\rho=\sqrt{x^2+(y+H)^2}.
\]
The bottom and right wall are fixed and impermeable, so homogeneous
Neumann conditions are imposed there.

This example requires a mesh map different from the vertical-column map
used for fixed sidewalls. Let $\eta_L$ be the elevation at the left contact
point. Its horizontal position and ALE velocity satisfy
\begin{equation}
    x_c=q(t)(H+\eta_L),
    \qquad
    w_c=\dot q(t)(H+\eta_L)+q(t)\dot\eta_L.
    \label{eq:rotating_contact_motion}
\end{equation}
Let $s_L=\eta_{h,x}$ denote the surface slope at the left endpoint.
The contact-point elevation then satisfies
\begin{equation}
    \dot\eta_L
    =
    \frac{(G_h)_L+\dot q(t)(H+\eta_L)s_L+(Q_\eta)_L}
         {1-q(t)s_L},
    \label{eq:rotating_contact_evolution}
\end{equation}
where $(Q_\eta)_L$ is zero when numerical dissipation is disabled.
For $P^1$, $s_L$ is the slope on the first surface element.
For $P^2$, it is the derivative of the surface trace polynomial
at the left endpoint.
For reference coordinates $s_i=i/N_x$ and vertical levels $r_j$ defined
in \eqref{SiRj}, the rotating-wall mesh vertices are placed at
\begin{equation}
    x_{i,j}=s_iL+r_j(1-s_i)x_c,
    \qquad
    y_{i,j}=-H+r_j(H+\eta_i).
    \label{eq:rotating_mesh_map}
\end{equation}
Here $\eta_i$ is the surface elevation at the $i$th top vertex.
The map places the left column on $x=q(t)(y+H)$ and the upper row
on the free surface, while keeping the bottom and right wall fixed. For $P^2$, top-edge midpoint heights are supplied by the
quadratic surface trace, and the remaining edge midpoints are placed
between their vertices. The contact point and interior mesh are updated
at each RK4 stage. A trial RK4 step is rejected if the minimum element Jacobian
falls below $0.02$ times its initial value. Rejection also occurs
if the mesh-displacement ratio exceeds $0.15$ or $1-q(t)s_L\leq0.1$.
A rejected step is retried with half the time step.

The parameters are set to
\begin{equation}
    L=4,\qquad H=1,\qquad g=1,\qquad
    A=0.025,\qquad \omega=2.
\end{equation}

\begin{figure}[htbp]
\centering
\includegraphics[width=0.48\textwidth]{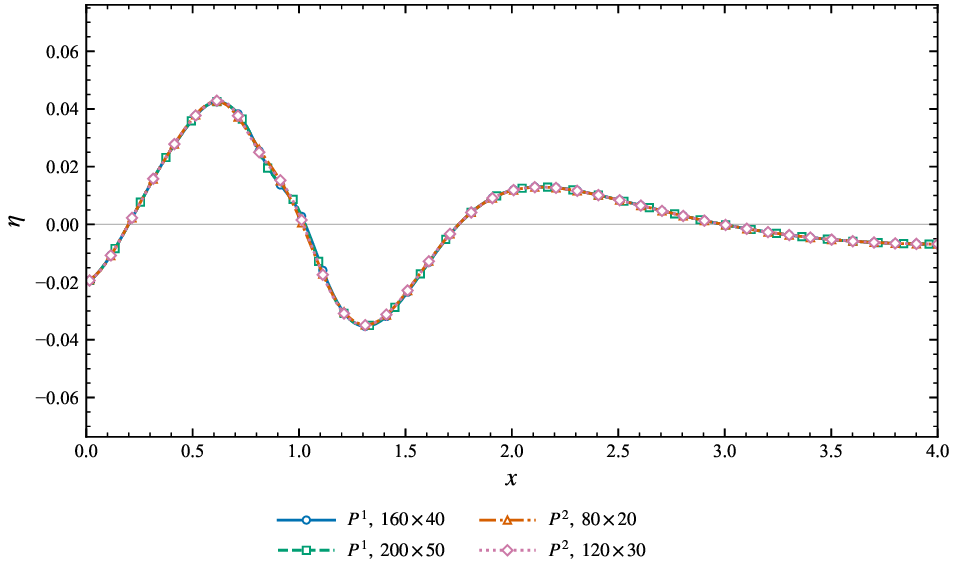}
\includegraphics[width=0.48\textwidth]{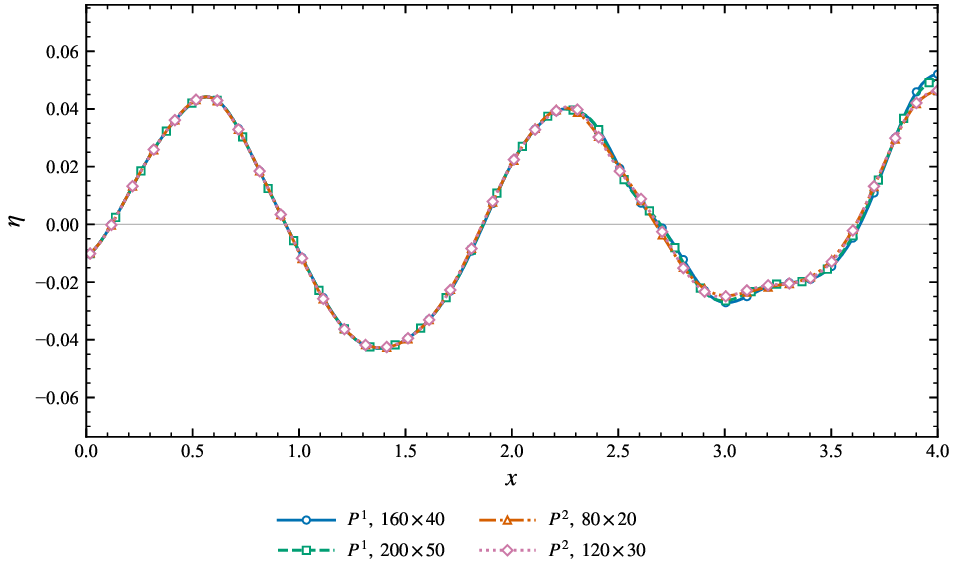}
\caption{Comparison of the free-surface elevations obtained with the $P^1$ and $P^2$ approximations on different mesh resolutions at  $t=7.5$ (top) and  $t=20$ (below).}
\label{fg:Non-eta}
\end{figure}

All results in this example include KO6 dissipation.
Figure~\ref{fg:Non-eta} compares the $P^1$ and $P^2$ free-surface
elevations at $t=7.5$ and $t=20$. The profiles agree closely, with
small local differences at the later time. Their closer agreement under refinement indicates reduced mesh
sensitivity. These profiles do not establish a quantitative
convergence rate.

\begin{figure}[htbp]
\centering
\begin{minipage}[t]{.62\textwidth}
\centering
\includegraphics[width=\linewidth]{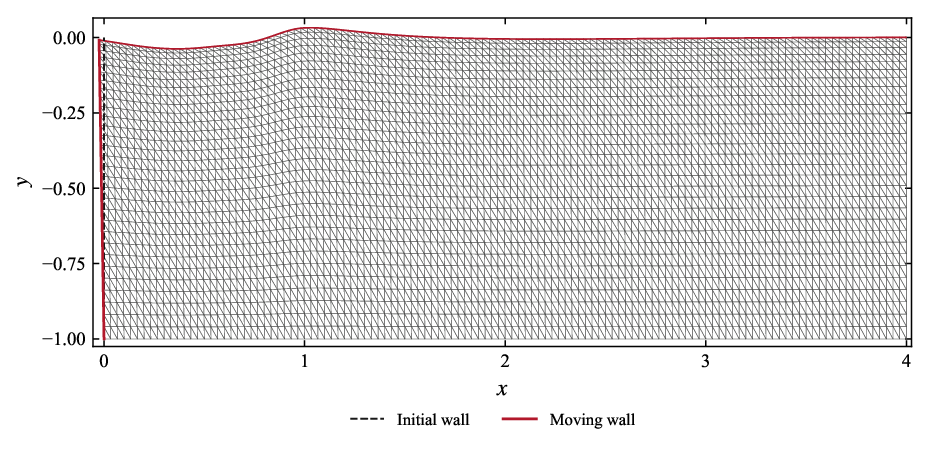}
\par\vspace{0.15em}
{\footnotesize (a) $t=5.5$, full mesh}
\end{minipage}
\hfill
\begin{minipage}[t]{.35\textwidth}
\centering
\includegraphics[width=\linewidth]{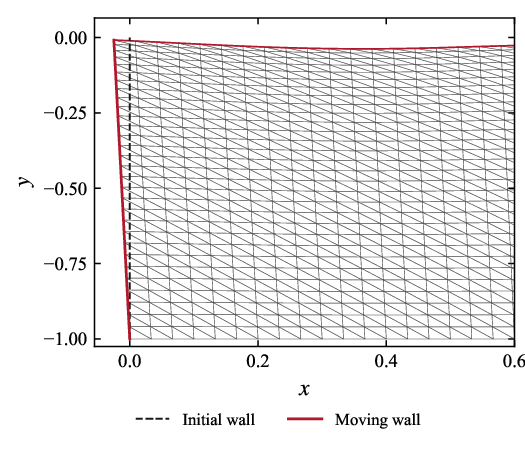}
\par\vspace{0.15em}
{\footnotesize (b) $t=5.5$, near the wavemaker}
\end{minipage}
\par\vspace{0.5em}
\begin{minipage}[t]{.62\textwidth}
\centering
\includegraphics[width=\linewidth]{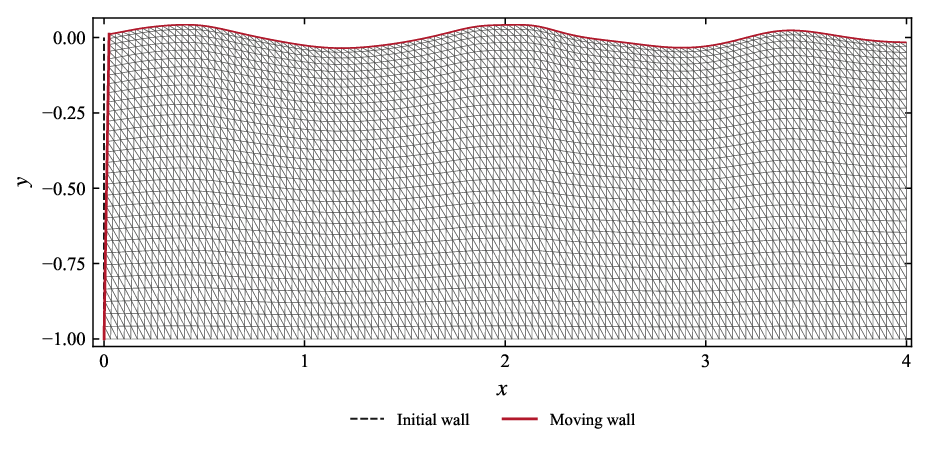}
\par\vspace{0.15em}
{\footnotesize (c) $t=16.5$, full mesh}
\end{minipage}
\hfill
\begin{minipage}[t]{.35\textwidth}
\centering
\includegraphics[width=\linewidth]{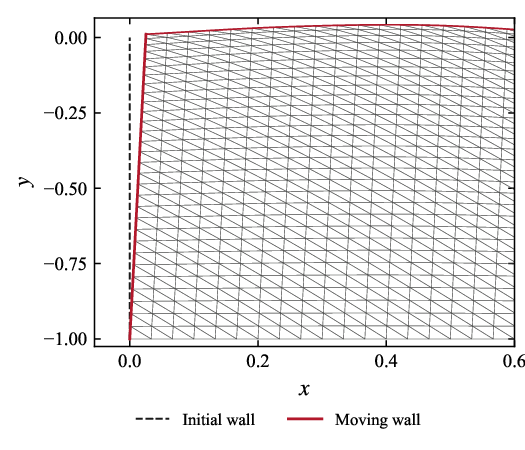}
\par\vspace{0.15em}
{\footnotesize (d) $t=16.5$, near the wavemaker}
\end{minipage}
\caption{ALE mesh connectivity for the $P^2$ approximation on the $120\times30$ grid at $t=5.5$ (top) and $t=16.5$ (bottom). Full meshes are shown on the left and wavemaker-region enlargements on the right.}
\label{fg:Non-potential}
\end{figure}

Figure~\ref{fg:Non-potential} shows the full mesh and local
enlargements at $t=5.5$ and $t=16.5$. The rotating-wall map is updated
at every RK4 stage; the bottom and right wall remain fixed while the
interior mesh deforms around the moving contact point.

\begin{figure}[htbp]
\centering
\includegraphics[width=0.7\textwidth]{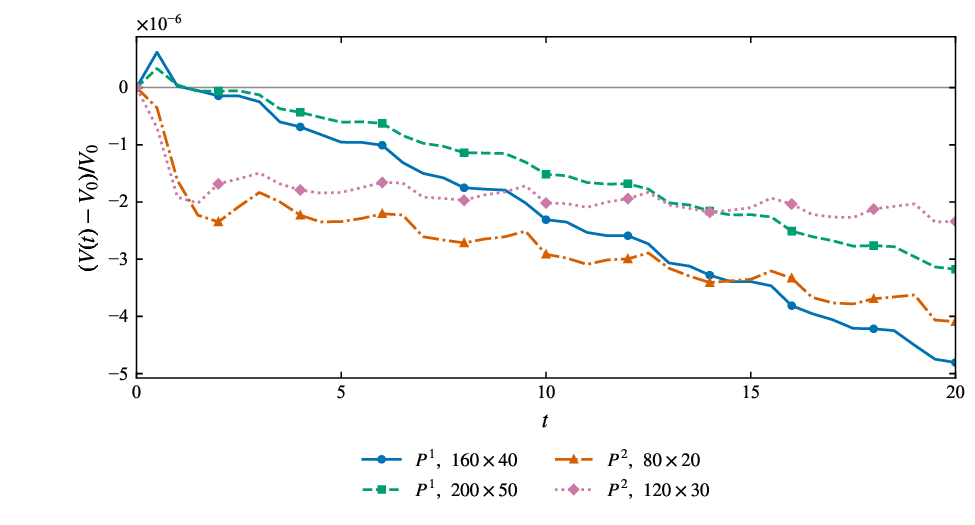}
\caption{Relative fluid-volume change for the rotating-wavemaker test with $P^1$ and $P^2$ approximations on four meshes. The initial volume is $V_0=4$.}
\label{fg:Non-volume}
\end{figure}

Figure~\ref{fg:Non-volume} plots the relative volume change
$(V(t)-V_0)/V_0$, where $V(t)$ includes the rotating-wall geometry and
$V_0=LH=4$. Diagnostics are recorded every $0.5$ time units for
$0.5\leq t\leq20$, with the value at $t=0$ taken from the initial
condition.

The semi-discrete formulation is volume preserving under the
compatibility assumptions stated above. After RK4 time discretization,
the relative volume change remains below $4.81\times10^{-6}$ for all
four discretizations, and refinement reduces the maximum drift for
each polynomial degree. Thus, the fully discrete RK4 method introduces
a small but nonzero volume error, while retaining volume conservation
to a very high numerical accuracy in this test.

\section{Conclusions}
\label{sec:conclusion}

We have developed an ALE continuous Galerkin finite element method
for two-dimensional fully nonlinear free-surface potential flow.
The method couples an isoparametric Laplace solver with ALE evolution
of the surface elevation and potential. Its main contribution is
to recover the Dirichlet--Neumann flux from the interior finite
element residual and connect its global balance to fluid-volume
conservation. Under the stated geometric and discretization
compatibility assumptions, this coupling conserves volume at the
semi-discrete level. A volume-tangent projection retains this
property when optional KO6 dissipation is introduced.
Numerical experiments demonstrate the expected
high-order behavior for linear waves, consistent results between
different polynomial degrees, and robust mesh deformation for curved
bottoms and a moving wavemaker.
The analysis remains conditional on the compatibility assumptions for
the ALE mesh motion. Future work will establish fully discrete
geometric volume identities and error estimates, develop
volume-preserving time integrators, and extend the method to nonlinear
stability analysis and three-dimensional flows.

\bibliographystyle{siamplain}
\bibliography{paper_final}

\end{document}